\documentclass[a4paper]{article}

\usepackage{amsthm,amsmath,amssymb}
\usepackage{mathrsfs}
\usepackage{graphicx}
\usepackage{float,verbatim}
\usepackage{threeparttable,booktabs}
\usepackage{underscore}
\usepackage{algorithm, algorithmicx}
\usepackage{subfigure}
\usepackage{diagbox}
\usepackage{multirow}
\usepackage{makecell}
\usepackage[noend]{algpseudocode}
\usepackage[nohead,margin=1.0in]{geometry}
\usepackage{color}
\usepackage{hyperref}
\hypersetup{hypertex=true,	colorlinks=true,	linkcolor=blue, anchorcolor=blue,        citecolor=blue}
\graphicspath{{./pic/}}
\usepackage{epstopdf}
\newtheorem{definition}{Definition}[section]

\newtheorem{proposition}{Proposition}[section]
\newtheorem{corollary}{Corollary}[section]
\newtheorem{theorem}{Theorem}[section]

\newtheorem{lemma}{Lemma}[section]
\newtheorem{example}{Example}[section]

\newtheorem{assumption}{Assumption}[section]

\algdef{SE}[DOWHILE]{Do}{doWhile}{\algorithmicdo}[1]{\algorithmicwhile\ #1}
\numberwithin{equation}{section}
\makeatletter
\@addtoreset{equation}{section}
\makeatother
\allowdisplaybreaks

\title{Point Source Identification Using Singularity Enriched Neural Networks\thanks{The work of B. Jin is supported by Hong Kong RGC General Research Fund (Project
14306423), and a start-up fund and Direct Grant for Research, both from The Chinese University of Hong Kong. The work of Z. Zhou is supported by Hong Kong
Research Grants Council (15303122) and an internal grant of Hong Kong Polytechnic University (Project ID: P0038888, Work
Programme: 1-ZVX3).}}
\author{Tianhao Hu\thanks{Department of Mathematics, The Chinese University of Hong Kong, Shatin, New Territories, Hong Kong, P.R. China  (\texttt{1155202953@link.cuhk.edu.hk}, \texttt{b.jin@cuhk.edu.hk})}\and Bangti Jin\footnotemark[2]\and Zhi Zhou\thanks{Department of Applied Mathematics, The Hong Kong Polytechnic University, Kowloon, Hong Kong, P.R. China (\texttt{zhizhou@polyu.edu.hk})}}
\date{}

\begin{document}
	\maketitle
  	
\begin{abstract}
The inverse problem of recovering point sources represents an important class of applied inverse problems. However, there is still a lack of neural network-based methods for point source identification, mainly due to the inherent solution singularity. In this work, we develop a novel algorithm to identify point sources, utilizing a neural network combined with a singularity enrichment technique. We employ the fundamental solution and neural networks to represent the singular and regular parts, respectively, and then minimize an empirical loss involving the intensities and locations of the unknown point sources, as well as the parameters of the neural network. Moreover, by combining the conditional stability argument of the inverse problem with the generalization error of the empirical loss, we conduct a rigorous error analysis of the algorithm. We demonstrate the effectiveness of the method with several challenging experiments.\\
\textbf{Key words:} point source identification, neural network, singularity enrichment, error estimate
\end{abstract}

\section{Introduction}\label{sec:intro}

Inverse source problems (ISP) represent an important class of applied inverse problems, and play an important role in science, engineering and medicine. One model setting is as follows. Let $\Omega\subset(-1,1)^d$ be an open bounded domain with a smooth  boundary $\partial\Omega$. Let $ \nu $ denote the unit outward normal vector to $\partial\Omega$, and $\partial_\nu u$ denote taking the outward normal derivative. Consider the standard Poisson equation \begin{equation}\label{problem}
    -\Delta u = F,\quad \mbox{in }\Omega.
\end{equation}
The concerned inverse problem is to determine the source   $F$ from the Cauchy data $(f,g)=(u|_{\partial\Omega}, \partial_\nu u|_{\partial\Omega})\in H^{\frac{3}{2}}(\partial\Omega)\times H^{\frac{1}{2}}(\partial\Omega)$. In practice, we can only access noisy data $(f^\delta,g^\delta)\in (L^2(\partial\Omega))^2$ with a noise level $\delta$, i.e. $\|f-f^\delta\|_{L^2(\partial\Omega)}\leq\delta$ and $\|g-g^\delta\|_{L^2(\partial\Omega)}\leq\delta$.

The ISP for the Poisson equation has been extensive studied (see the monographs \cite{anikonov2013inverse, isakov1990inverse}). With only one Cauchy data pair, it is impossible to uniquely recover a general source $F$: Indeed, adding any smooth function $v$ with a compact support $\overline{\mathrm{supp}(v)}\subset \Omega$ to $u$ does not change the Cauchy data. Thus, suitable prior knowledge (i.e., physical nature) on the source $F$ is required in order to restore unique identifiability. It is known that one Cauchy data pair can uniquely determine the unknown source $F$, e.g., one factor of a product source $F$ (i.e. $F$ is separable and one factor is known) \cite{kriegsmann1988source,el1998some}, or the angular factor of a spherically layered source \cite{el1998some}, or $F = \chi_D\rho(x)$ or $F=\nabla \cdot (\rho(x)\chi_D(x)a)$ (with $\chi_D$ being the characteristic function of a polygon $D\subset\Omega$, and $a$ a nonzero constant vector) \cite{ikehata1999reconstruction}.

In this work, we investigate the case that $F$ is a linear combination of point sources:
\begin{equation}\label{source}
    F = \sum_{j=1}^Mc_j\delta_{\mathbf{x}_j},\quad \mathbf{x}_j\in\Omega,\quad j=1,2,\cdots,M,
\end{equation}
where $\delta_{\mathbf{x}_j}$ is the Dirac delta function concentrated at $\mathbf{x}_j\in\Omega$, defined by $\int_\Omega \delta_{\mathbf{x}_j}v {\rm d}\mathbf{x} = v(\mathbf{x}_j)$ for all $v\in C(\overline{\Omega})$, and $c_j\neq0$ is the intensity of the point source at the point $\mathbf{x}_j$. The unknown parameters include the number $M$ of point sources, and the locations and intensities $\{(\mathbf{x}_j, c_j)\}_{j=1}^M$.
It arises in several practical applications, e.g., in electroencephalography, the Dirac function can describe active sources in the brain \cite{MICHEL20042195, michel2004128}.
The specific inverse problem satisfies unique identifiability, and H\"{o}lder stability \cite{ElBadiaElHajj:2012,Badia2013}.

Numerically, El Badia et al \cite{el2000inverse} developed a direct method to identify the number and locations / intensities for the Poisson problem in 2D domains.
From the Cauchy data $(f^\delta,g^\delta)$, we can compute the value $\sum_{j=1}^Mc_jw(\mathbf{x}_j)$ for any harmonic $w$. By suitably choosing $w$, we can determine the number $M$, and locations and densities of the point sources. It works well for both 2D and 3D cases, for which $w$ is constructed from complex analytic functions. It is more involved in the high-dimensional case  \cite{chafik2000some}: both constructing the harmonic test function $w$ and  computing the reciprocity gap functional (cf. \eqref{eq:R-formula}) are nontrivial. We aim at developing a novel neural method for the task using neural networks (NNs). The key challenge lies in the low regularity of point sources, and the low regularity of the associated PDE solution. Thus, standard neural solvers are ineffective in resolving the solution $u(\mathbf{x})$.
To overcome the challenge, we draw on the analytic insight that  $u(\mathbf{x})$ can be split into a singular part (represented by fundamental solutions), and a regular part that can be effectively approximated by NNs.

In this work, we develop a novel numerical method for recovering point sources using NNs. The proposed method consists of two steps. First we detect the number $M$ of point sources as \cite{el2000inverse}, and then learn simultaneously the smooth part (approximated using NNs) and parameters of the unknown point sources (i.e., the intensities and locations) using a least-squares type empirical loss. In sum, we make the following contributions.
First, we develop an easy-to-implement algorithm for recovering point sources. Second, we provide error bounds for the algorithm. This is achieved by suitably combining conditional stability with \textit{a priori} estimates on the empirical loss, using techniques from statistical learning theory \cite{AnthonyBartlett:1999}. We bound the Hausdorff distance between the approximate and exact locations, and errors of the intensities and the regular part, explicitly in terms of the noise level, neural network architecture and the numbers of sampling points in the domain and on the boundary.
Third, we present several numerical experiments, including partial Cauchy data and Helmholtz problem, to illustrate its flexibility and accuracy.

Now let us put the work into the context of neural solvers for ISPs. Solving PDE inverse problems using neural networks has received much recent attention in either supervised or unsupervised manners.
Nonetheless, neural solvers for ISPs are still very limited. Zhang and Liu \cite{Zhang2023} employed the deep Galerkin method \cite{sirignano2018dgm} to recover the source term in an elliptic problem using the measurement of the solution in a subdomain, and established the convergence of the method for analytic sources; see also \cite{ZhangLiLiu:2023} for the parabolic case. Li et al \cite{LiLiangWang:2023} investigated the use of DNNs for the stochastic inverse source problem. In contrast, this work focuses on identifying point sources which is strongly nonlinear but enjoys H\"{o}lder stability, which enables deriving \textit{a priori} error estimates for the proposed neural solver. Du et al \cite{DuLiSun:2023} developed a novel divide-and-conquer algorithm for recovering point sources for the Helmholtz equation. In contrast, this work contributes one novel unsupervised neural solver for identifying point sources and provides rigorous mathematical guarantee. 

The rest of the paper is organized as follows. We develop the method in Section \ref{sec:algorithm}. In Section \ref{sec:error}, we provide an analysis of the method. In Section \ref{sec:experiment}, we present numerical experiments to illustrate the performance of the method, and give concluding remarks in Section \ref{sec:concl}. In the appendices, we collect several technical estimates for analyzing the empirical loss.
Throughout, we use the notation $C$ to denote a constant which depends on its argument and whose value may differ at different occurrences.

\section{The proposed method}\label{sec:algorithm}
Now we describe the proposed method, termed as singularity enriched neural network (SENN), for identifying point sources. The procedure consists of two steps: (i) detect the number $M$ of point sources and (ii) estimate the locations and intensities with singularity enrichment. Throughout, for a harmonic function $w$, we denote
by $\mathcal{R}(w;f^\delta,g^\delta)$ the following functional
\begin{equation}\label{eq:R-formula}
    \mathcal{R}(w;f^\delta,g^\delta)\triangleq\int_{\partial\Omega}-g^\delta w+f^\delta\partial_\nu w{\rm d}\mathbf{x},
\end{equation}
and often abbreviate it to $\mathcal{R}(w)$. It is commonly known as the reciprocity gap functional.  It will play a key role in establishing the stability estimate.

\subsection{Detection of the number of sources}
First we recall the detection of the number $M$ of point sources. We only describe the case $d=2$, following \cite{el2000inverse}.  Note that the real or imaginary part of a complex polynomial is harmonic. By setting  $w_m(\mathbf{x})=w_m(x,y)=z^m=(x+{\rm i}y)^m,m\in\mathbb{N}$, we obtain
$$\mathcal{R}(w_m)=\sum_{j=1}^Mc_j(x_j+{\rm i}y_j)^m.$$
Given an upper bound $\overline{M}\geq M$, we can find  $M$ using the next result \cite[Lemma 2]{el2000inverse}:

\begin{proposition}
    \label{prop:number}
Let the Hankel matrix $A=(a_{ij})_{\overline{M}\times\overline{M}}$, with $a_{ij}=\mathcal{R}(w_{i+j-2})\in\mathbb{C}$. Then ${\rm rank}(A)=M$ and the first $M$ columns of $A$ are linearly independent.
\end{proposition}

Proposition \ref{prop:number} suggests that we may start with the first column of $A$, work forward one column each time, and stop the process when the columns become linearly dependent. The stopping index gives $M$. In practice, due to the presence of quadrature error in evaluating the entries $a_{ij}$, one should enforce the stopping criterion inexactly. Let $\mu_k$ be the $k$-th column of $A$, $A_k=(\mu_1,\mu_2,\cdots,\mu_k)$, and let $\sigma_s(A)$ be the smallest singular value of $A$. If $\sigma_s(A_k)<\epsilon_{\rm tol}$, a given tolerance, then we claim that the columns $\mu_1,\ldots,\mu_k$ are linearly dependent. The whole procedure is shown in Algorithm \ref{algorithm:num}.

\begin{algorithm}[H]
\caption{Detection of the number $M$ of sources in the 2d case}\label{algorithm:num}
\begin{algorithmic}[1]
\State Set $k=0$, upper bound $\overline{M}$, and the tolerance $\epsilon_{\rm tol}$.
\State Compute $\mathcal{R}(w_m)$ using \eqref{eq:R-formula} for $m=0,1,\cdots,\overline{M}-1$.
\State Set $A=(\mathcal{R}(w_0), \mathcal{R}(w_1), \cdots, \mathcal{R}(w_{\overline{M}-1}))^\top$.
\Do
    \State $m=m+1$, $k=k+1$.
    \State Compute $\mathcal{R}(w_m)$ using \eqref{eq:R-formula}. Let $\mu=(\mathcal{R}(w_k), \mathcal{R}(w_{k+1}), \cdots, \mathcal{R}(w_{k+\overline{M}-1}))^\top$ and set $A=(A,\mu)$.
    \State Compute $\sigma_s(A)$.
\doWhile{$\sigma_s(A)\geq\epsilon_{\rm tol}$.}
\State Output $M=k$.
	\end{algorithmic}
\end{algorithm}

For the stability of the method, we define the Vandermonde matrix $B\in\mathbb{C}^{\overline{M}\times M}$ by
\begin{equation*}
    B=\begin{pmatrix}
            1&1&\cdots&1\\
            x_1+{\rm i}y_1&x_2+{\rm i}y_2&\cdots&x_M+{\rm i}y_M\\
            \vdots&\vdots&\ddots&\vdots\\
            (x_1+{\rm i}y_1)^{\overline{M}-1}&(x_2+{\rm i}y_2)^{\overline{M}-1}&\cdots&(x_M+{\rm i}y_M)^{\overline{M}-1}\\
        \end{pmatrix}
        :=\begin{pmatrix}
            \\
            \text{\large{$(B_M)$}}_{\text{\tiny{$M\times M$}}}\\
            \\
            \text{\large{$(B_{\widetilde{M}})$}}_{\text{\tiny{$(\overline{M}-M)\times M$}}}
            \\
            \\
        \end{pmatrix},
    \end{equation*}
and let  $D={\rm diag}(x_1+{\rm i}y_1, x_2+{\rm i}y_2, \cdots, x_M+{\rm i}y_M)$. The $k$-th column of $A$ is $\mu_k=BD^{k-1}\mathbf{c}\in\mathbb{C}^M$, and $A_M=(\mu_1,\mu_2\cdots,\mu_M)=B(\mathbf{c}, D\mathbf{c},\cdots,D^{M-1}\mathbf{c})$. Let $C=(\mathbf{c}, D\mathbf{c},\cdots,D^{M-1}\mathbf{c})$. Then we have
    $ A_M=BC=\begin{pmatrix}
        B_MC\\
        B_{\widetilde{M}}C
    \end{pmatrix}. $
Since $B_MC$ is a submatrix of $A_M$, we have $\sigma_s(A_M)\geq\sigma_s(B_MC)>0$, which is independent of $\overline{M}$. Thus, in theory, we may pick the tolerance $\epsilon_{\rm tol}\leq\sigma_s(B_MC)$ in Algorithm \ref{algorithm:num}, and in practice, fix it at about the noise level (including quadrature error etc).

\subsection{Singularity enrichment}
Once the number $M$ is estimated, we apply deep learning to learn the solution $u$ and the locations and densities. Due to the presence of point sources, $u$ is nonsmooth, and standard neural PDE solvers are ineffective \cite{HuJinZhou:2022}. We employ singularity enrichment, which splits out the singular part, in order to resolve the challenge. This idea has been widely used in the context of FEM \cite{Fix:1973, cai2001finite,CaiShin:2001, Fries:2010}, and also neural PDE solvers \cite{HuJinZhou:2022,hu2023solving}. We extend the idea to recover point sources.

We employ the fundamental solution $\Phi(\mathbf{x})$ to the Laplace equation \cite[p. 22]{evans10}:
\begin{equation}\label{fundamental-s}
    \Phi(\mathbf{x})=c_d\left\{
	\begin{aligned}
		-\log|\mathbf{x}|,&\quad d=2,\\
		|\mathbf{x}|^{2-d},&\quad d\geq3,\\
	\end{aligned}\right.
\end{equation}
with $c_2=\frac{1}{2\pi}$ and $c_d=\frac{1}{d(d-2)\alpha(d)}$ for $d\geq3$,
with $ \alpha(d) =\pi^\frac{d}{2}/\Gamma(\frac{d}{2}+1)$. 
By definition, for any $\mathbf{x}_0\in\mathbb{R}^d$, $ \Phi(\mathbf{x}) $ satisfies
$-\Delta \Phi(\mathbf{x-x}_0)=\delta_{\mathbf{x}_0}$ in  $\mathbb{R}^d\setminus\{\mathbf{x}_0\}$.
We abbreviate $\Phi(\mathbf{x}-\mathbf{x}_0)$ as $\Phi_{\mathbf{x}_0}$. Using $\Phi(\mathbf{x})$, we can split the solution $u$ to \eqref{problem} as
\begin{equation}\label{split}
    	u=\sum_{j=1}^Mc_j\Phi_{\mathbf{x}_j}+v.
\end{equation}
The term $\sum_{j=1}^Mc_j\Phi_{\mathbf{x}_j}$ captures the leading singularity of the solution $u$ due to the point sources $\sum_{j=1}^Mc_j\delta_{\mathbf{x}_j}$ and $ v $ is the regular part. The regular part $v$ satisfies
\begin{equation}\label{problem-v}
    \left\{
    \begin{aligned}
    -\Delta v&=0,&& \mbox{in }\Omega,\\
    v&=f-F(\mathbf{c},X),&& \mbox{on }\partial\Omega,\\
    \partial_\nu v&=g-\partial_\nu F(\mathbf{c},X),&& \mbox{on }\partial\Omega,
    \end{aligned}
    \right.\quad \mbox{with } F(\mathbf{c},X) =\sum_{j=1}^Mc_j\Phi_{\mathbf{x}_j}.
\end{equation}
In \eqref{split}, we have used the information about the unknown point source.
This reformulation motivates the following loss for noisy data:
\begin{equation*}
    L^\delta(v,\mathbf{c},X)=\|\Delta v\|_{L^2(\Omega)}^2 +\sigma_d\big\|v-f^\delta+F(\mathbf{c},X)\big\|_{L^2(\partial\Omega)}^2+\sigma_n\big\|\partial_\nu v-g^\delta+\partial_\nu F(\mathbf{c},X)\big\|_{L^2(\partial\Omega)}^2,
\end{equation*}
where $\boldsymbol\sigma=(\sigma_d,\sigma_n)^\top\in\mathbb{R}^2_+$ is the weight vector, $\mathbf{c}=(c_1,c_2,\cdots,c_M)^\top$ the vector of intensities, and $X=(\mathbf{x}_1, \mathbf{x}_2,\cdots, \mathbf{x}_M)$ the matrix of locations.
We denote the loss for the exact data $(f,g)$ by $L(v,\mathbf{c},X)$, and its minimizer by  $(v^*, \mathbf{c}^*, X^*)$. Note that $L(v^*,\mathbf{c}^*,X^*)=0$.

\subsection{Reconstruction via SENN}
To approximate the regular part $v$, we employ standard fully connected feedforward NNs. These are functions $ f_\theta: \mathbb{R}^d\rightarrow \mathbb{R}$, with the parameter $\theta$. It is defined recursively by
$ \mathbf{x}_0=\mathbf{x}$, $\mathbf{x}_\ell=\varrho(\mathbf{A}_\ell\mathbf{x}_{\ell-1}+\mathbf{b}_\ell)$, $\ell=1,2,\cdots,L-1$,	and $f_\theta(\mathbf{x})=\mathbf{A}_L\mathbf{x}_{L-1}+\mathbf{b}_L$,
where $\mathbf{A}_\ell\in\mathbb{R}^{n_\ell\times n_{\ell-1}}$, $\mathbf{b}_\ell\in\mathbb{R}^{n_\ell},$ $\ell=1,2,\cdots,L$, with $n_0=d$ and $n_L=1$. $\varrho:\mathbb{R}\to \mathbb{R}$ is a nonlinear function, applied componentwise to a vector, and we employ the hyperbolic tangent $\varrho (x)=\frac{{\rm e}^x-{\rm e}^{-x}}{{\rm e}^x+{\rm e}^{-x}}$. The integers $L$ and $W:=\max_{\ell=0,\ldots,L} n_\ell$ are the depth and width. The parameters $\mathbf{A}_\ell$ and $\mathbf{b}_\ell\,\,(\ell=1,2,\cdots,L)$ are trainable and collected in a vector $\theta$. We denote by $\mathcal{N}_\varrho(L,N_\theta,B_\theta)$ the set of NN functions of depth $L$, with the total number of nonzero parameters $N_\theta$, and each entry bounded by $B_\theta$, using the activation function $\varrho$:
\begin{equation}\label{eqn:nn-set}
  \mathcal{N}_\varrho(L,N_\theta,B_\theta)   = \{v_\theta: v_\theta ~ \mbox{has a depth } L,\, |\theta|_0\leq N_\theta, \, |\theta|_{\ell^\infty}\leq B_\theta\},
\end{equation}
where $|\cdot|_{\ell^0}$ and $|\cdot|_{\ell^\infty}$ denote the number of nonzero entries in and the maximum norm of a vector, respectively. Below we also use the shorthand notation $\mathcal{A}=\mathcal{N}_\varrho(L,N_\theta,B_\theta)$.

We employ an element $v_\theta\in \mathcal{A}$ to approximate the regular part $v$, and  discretize relevant integrals using quadrature, e.g., Monte Carlo. Let $U(D)$ be the uniform distribution over a set $D$, and let $|D|$ denote its Lebesgue measure. Then we can rewrite the loss $\widehat{L}^\delta(v_\theta,\mathbf{c},X)$ as
\begin{align*}
L^\delta(v_\theta,\mathbf{c},X)=&|\Omega|\mathbb{E}_{ U(\Omega)}[(\Delta v_\theta(Z))^2] +\sigma_d|\partial\Omega|\mathbb{E}_{U(\partial\Omega)}\big[\big((v_\theta-f^\delta+F(\mathbf{c},X))(Y)\big)^2\big]\\
&+\sigma_n|\partial\Omega|\mathbb{E}_{U(\partial\Omega)}\big[\big((\partial_\nu v_\theta-g^\delta+\partial_\nu F(\mathbf{c},X))(Y)\big)^2\big].
\end{align*}
where $\mathbb{E}_\nu$ denotes taking expectation with respect to $\nu$. Let $\{X_{i_1}\}_{i_1=1}^{N_r}$ and $\{Y_{i_2}\}_{i_2=1}^{N_b}$ be identically and independently distributed (i.i.d.) sampling points drawn from $U(\Omega)$ and $U(\partial\Omega)$, respectively, i.e., $\{X_{i_1}\}_{i_1=1}^{N_r}\sim U(\Omega)$ and $\{Y_{i_2}\}_{i_2=1}^{N_b}\sim U(\partial\Omega)$.
Then the empirical loss $\widehat{L}^\delta$ is given by
\begin{align}
\widehat{L}^\delta(v_\theta,\mathbf{c},X)=&\dfrac{|\Omega|}{N_r}\sum_{i_1=1}^{N_r}|\Delta v_\theta(X_{i_1})|^2+\sigma_d\dfrac{|\partial\Omega|}{N_{b}}\sum_{i_2=1}^{N_b}\big(v_\theta(Y_{i_2})-f^\delta(Y_{i_2})+F(\mathbf{c},X)(Y_{i_2})\big)^2\nonumber\\
&+\sigma_n\dfrac{|\partial\Omega|}{N_{b}}\sum_{i_2=1}^{N_b}\big(\partial_\nu v_\theta(Y_{i_2})-g^\delta(Y_{i_2})+\partial_\nu F(\mathbf{c},X)(Y_{i_2})\big)^2.
\label{eqn:loss-emp}
\end{align}

The resulting training problem reads
\begin{equation}\label{eqn:min-emp}
    (\widehat{\theta}^*,\widehat{\mathbf{c}}^*,\widehat{X}^*) = \mathop{{\rm argmin}}\limits_{\theta,\mathbf{c},X}\widehat{L}^\delta(v_\theta,\mathbf{c},X),
\end{equation}
with $v_{\widehat{\theta}^*}\in\mathcal{A}$ being the NN approximation. Since $(\theta,\mathbf{c},X)$
is finite-dimensional and lives on a compact set (due to the $\ell^\infty$ bound $|\theta|_{\infty}\leq B_\theta$), and we further assume $|\mathbf{c}|_\infty\leq B_\mathbf{c}$, for smooth $\varrho$, the loss $\widehat {L}^\delta(v_\theta,\mathbf{c},X)$ is continuous in $(\theta,\mathbf{c},X)$, which directly implies the existence of a global minimizer. Problem \eqref{eqn:min-emp} is nonconvex, due to the nonlinearity of the NN $v_\theta$ in $\theta$. Hence, it may be challenging to find a global minimizer of $\widehat{L}^\delta(v_\theta,\mathbf{c},X)$. In practice, gradient type algorithms, e.g.,  Adam \cite{KingmaBa:2015} or L-BFGS \cite{ByrdLu:1995}, seem to work fairly well \cite{Karniadakis:2021nature}. The gradient of $\widehat{L}(u_\theta)$ with respect to the NN parameters $\theta$  and the NN $u_\theta(\mathbf{x})$ with respect to the input $\mathbf{x}$ can be computed via automatic differentiation. Thus, the method is easy to implement.

\section{Error analysis}\label{sec:error}
Now we give an error analysis of the method in Section \ref{sec:algorithm} in order to provide theoretical guarantee. The analysis is challenging due to the nonlinearity of the ISP. Our analysis combines conditional stability of the ISP with the consistency error analysis for the empirical loss. For neural solvers for PDE direct problems, the analysis of the empirical loss has been studied (see, e.g., \cite{JiaoLai:2022cicp,LuChenLu:2021,HuJinZhou:2022}).

\subsection{Conditional stability}\label{ssec:stability}
First we analyze the stability of the loss $L(v_\theta, \mathbf{c}, X)$. Throughout, let $\mathbf{x}=(x_1,x_2,\cdots,x_d)^\top \in \Omega\subset\mathbb{R}^d$ and ${\rm dist}(\mathbf{x}, \partial\Omega)={\inf_{\mathbf{y}\in\partial\Omega}}{\rm dist}(\mathbf{x},\mathbf{y})$ be the Euclidean distance between $\mathbf{x}$ and $\partial\Omega$. Let $\gamma={\min_{1\leq j\leq M}}{\rm dist}(\mathbf{x}_j, \partial\Omega)>0$ be the minimal distance of the set of points $\{\mathbf{x}_j\}_{j=1}^M$ to $\partial\Omega$, and let
\begin{equation}
   \Omega_\gamma= \{{\mathbf{x}\in\Omega:{\rm dist}(\mathbf{x},\partial\Omega)}\geq\gamma\}.
\end{equation}
Then $\{\mathbf{x}_j\}_{j=1}^M\subset \Omega_\gamma$. Let $\beta={\rm diam}(\Omega)-\gamma$, with $\mathrm{diam}(\Omega)$ being the diameter of  $\Omega$.
Further, we assume that $0<b_\mathbf{c}\leq |c_j^*|\leq B_\mathbf{c}$, $j=1,\ldots,M$. Let $(\widehat{\theta}^*,\widehat{\mathbf{c}}^*,\widehat{X}^*)$ be a minimizer of $\widehat{L}^\delta(v_\theta,\mathbf{c},X)$ over the set $\mathbb{A} = \mathcal{A}\times I_{\mathbf{c}}^M \times \Omega_\gamma^M$ with $I_{\mathbf{c}} = [-B_\mathbf{c},-b_\mathbf{c}]\cup[b_\mathbf{c},B_\mathbf{c}]$. 

First we bound the error of the empirical risk minimizer $\widehat{X}^*$ relative to the exact one $X^*$ using the population loss value $L(v_{\widehat{\theta}^*},\widehat{\mathbf{c}}^*,\widehat{X}^*)$. The proof is inspired by the conditional stability analysis of the ISP in \cite{ElBadiaElHajj:2012,Badia2013}. The key is to multiply with equation \eqref{problem} some carefully constructed test functions and to exploit the reciprocity gap functional $\mathcal{R}(w;f,g)$. Then using analytical properties of Vandemonde type matrices, we can obtain conditional stability for identifying point sources from Cauchy data. The construction relies on the harmonicity of complex analytic functions, which exist only in the two-dimensional complex plane, and we project the points in the general $d$-dimensional domain into two-dimensional complex plane. More specifically, for one configuration $X=(\mathbf{x}_1, \mathbf{x}_2, \cdots, \mathbf{x}_M)$, with $\mathbf{x}_j=(x_{1,j},x_{2,j},\cdots,x_{d,j})^\top$, we denote by $P_{k}(\mathbf{x}_j)=x_{k,j}+{\rm i}x_{k+1,j}$, $k=1,\ldots,d-1$, to be its projection onto $x_kOx_{k+1}$-complex plane and let $\mathbf{P}_{k}(X)=(P_{k}(\mathbf{x}_1),P_{k}(\mathbf{x}_2),\cdots,P_{k}(\mathbf{x}_M))$. Then we define the separability coefficient $\rho_k(X)$ of the projected sources $\mathbf{P}_k(X)$ and the distance $\rho(X)$ by
\begin{equation}\label{def:sepa}
    \rho_k(X)=\min_{1\leq i<j\leq M}{\rm dist}(P_{k}(\mathbf{x}_i),P_{k}(\mathbf{x}_j))\quad\mbox{and}\quad\rho(X)=\min_{1\leq k<d}\rho_k(X).
\end{equation}
For two point configurations $\widehat{X}^*=(\widehat{\mathbf{x}}_1^*, \widehat{\mathbf{x}}_2^*,\cdots, \widehat{\mathbf{x}}_{\widehat{M}}^*)$ and $X^*=\left(\mathbf{x}_1^*, \mathbf{x}_2^*,\cdots, \mathbf{x}_{M}^*\right)$, we make the following separability assumption for the subsequent analysis.
\begin{assumption}\label{assump:separation}
There exists a $\rho_*>0$ such that
$ \min\{\rho(\widehat{X}^*), \rho(X^*)\}\geq\rho_*>0$.
\end{assumption}

We define the Hausdorff distance $d_H(\widehat{X}^*,X^*)$ by
\begin{equation*}
    d_H(\widehat{X}^*,X^*)=\max\Big\{\max_{1\leq i\leq M}\min_{1\leq j\leq\widehat{M}}{\rm dist}\left(\mathbf{x}_i^*, \widehat{\mathbf{x}}_j^*\right), \max_{1\leq j\leq \widehat{M}}\min_{1\leq i\leq M}{\rm dist}\left( \widehat{\mathbf{x}}_j^*,\mathbf{x}_i^*\right)\Big\}.
\end{equation*}

We have the following stability result in Hausdorff distance.
\begin{theorem}\label{thm:disbound}
Let $X^*$ and $\widehat{X}^*$ be minimizers of the losses ${L}$ and $\widehat{L}^\delta$, respectively. Then under Assumption \ref{assump:separation}, there holds
    \begin{equation*}
        d_H(\widehat{X}^*, X^*)\leq {C}(\Omega,\beta,M,\widehat{M},\boldsymbol{\sigma},b_\mathbf{c},\rho_*,d)L(v_{\widehat{\theta}^*},\widehat{\mathbf{c}}^*,\widehat{X}^*)^{\frac{1}{2\max(M,\widehat{M})}}.
    \end{equation*}
\end{theorem}
\begin{proof}
Let $u_{\widehat{\theta}^*}=v_{\widehat{\theta}^*}+\sum_{j=1}^{\widehat{M}}\widehat{c}_j^*\Phi_{\widehat{\mathbf{x}}_j^*}\left(\mathbf{x}\right)$. Then  $u^*$ and $u_{\widehat{\theta}^*}$ satisfy
\begin{equation}\label{equation-num-true}
    \left\{ {\begin{aligned}
     -\Delta u^* &= \sum_{j=1}^Mc_j^*\delta_{\mathbf{x}_j^*}, &&\mbox{in }\Omega,\\
     u^*&=f, &&\mbox{on }\partial\Omega,\\
     \partial_\nu u^*&=g, &&\mbox{on }\partial\Omega,
    \end{aligned}} \right.\quad\text{and}\quad
    \left\{ {\begin{aligned}
     -\Delta u_{\widehat{\theta}^*} &=-\Delta v_{\widehat{\theta}^*}+\sum_{j=1}^{\widehat{M}}\widehat{c}_j^*\delta_{\widehat{\mathbf{x}}_j^*}, &&\mbox{in }\Omega,\\
     u_{\widehat{\theta}^*}&=\widehat{f}, &&\mbox{on }\partial\Omega,\\
     \partial_\nu u_{\widehat{\theta}^*}&=\widehat{g}, &&\mbox{on }\partial\Omega,
    \end{aligned}} \right.
    \end{equation}
respectively, with $\widehat f$ and $\widehat g$ given by
\begin{equation}\label{eqn:hatfg}
    \widehat{f}= v_{\widehat{\theta}^*}+\sum_{j=1}^{\widehat{M}}\widehat{c}_j^*\Phi_{\widehat{\mathbf{x}}_j^*}\quad \mbox{and}\quad \widehat{g}\triangleq\partial_\nu v_{\widehat{\theta}^*}+\sum_{j=1}^{\widehat{M}}\widehat{c}_j^*\partial_\nu\Phi_{\widehat{\mathbf{x}}_j^*}.
\end{equation}
Now consider the $d$-dimensional complex harmonic functions:
\begin{equation}
\Psi_n(\mathbf{x})=\prod_{j=1,j\neq n}^{M}\left(P_k(\mathbf{x})-P_{k}(\mathbf{x}_j^*)\right)\prod_{j=1}^{\widehat{M}}(P_k(\mathbf{x})-P_{k}(\widehat{\mathbf{x}}_j^*)),\quad n=1,2,\cdots,M.
\end{equation}
Setting $w=\Psi_n $ in the reciprocity gap functional $\mathcal{R}(w;f,g)$ (cf. \eqref{eq:R-formula}) gives
\begin{equation*}
    \mathcal{R}(\Psi_n;f,g)=\sum_{j=1}^Mc_j^*\Psi_n(\mathbf{x}_{j}^*)=c_n^*\Psi_n(\mathbf{x}_{n}^*).
\end{equation*}
Next multiplying $-\Delta u_{\widehat{\theta}^*}=-\Delta v_{\widehat{\theta}^*}+\sum_{j=1}^{\widehat{M}}\widehat{c}_j^*\delta_{\widehat{\mathbf{x}}_j^*}$ with $\Psi_n$ and then integrating by parts twice yield
\begin{equation*}
    \mathcal{R}(\Psi_n;\widehat{f},\widehat{g})=\sum_{j=1}^{\widehat{M}}\widehat{c}_j^*\Psi_n(\widehat{\mathbf{x}}_{j}^*)-\int_\Omega\Delta v_{\widehat{\theta}^*}\Psi_n{\rm d}\mathbf{x}=-\int_\Omega\Delta v_{\widehat{\theta}^*}\Psi_n{\rm d}\mathbf{x}.
    \end{equation*}
    Subtracting the preceding two identities yields that for $n=1,2,\cdots,M$,
    $$c_n^*\Psi_n(\mathbf{x}_{n}^*)=\mathcal{R}(\Psi_n;f-\widehat{f},g-\widehat{g})-\int_\Omega\Delta v_{\widehat{\theta}^*}\Psi_n{\rm d}\mathbf{x},$$
    which, upon substituting the definition of  $\Psi_n$, is equivalent to
\begin{equation*}
        c_n^*\prod_{j=1,j\neq n}^{M}\!\!\!\!\!\left(P_k(\mathbf{x}_n^*)-P_{k}(\mathbf{x}_j^*)\right)\prod_{j=1}^{\widehat{M}}\left(P_k(\mathbf{x}_n^*)-P_{k}(\widehat{\mathbf{x}}_j^*)\right)=\mathcal{R}(\Psi_n;f-\widehat{f},g-\widehat{g})-\int_\Omega\Delta v_{\widehat{\theta}^*}\Psi_n{\rm d}\mathbf{x}.
    \end{equation*}
    Using the separability coefficient $\rho_k(X)$ defined in \eqref{def:sepa} and Assumption \ref{assump:separation}, we have
\begin{align*}
        &\bigg|c_n^*\prod_{j=1,j\neq n}^{M}\left(P_k(\mathbf{x}_n^*)-P_{k}(\mathbf{x}_j^*)\right)\prod_{j=1}^{\widehat{M}}\left(P_k(\mathbf{x}_n^*)-P_{k}(\widehat{\mathbf{x}}_j^*)\right)\bigg|\\
        \geq&b_\mathbf{c}\Big(\min_{1\leq j\leq \widehat{M}}{\rm dist}(P_k(\mathbf{x}_n^*),P_{k}(\widehat{\mathbf{x}}_j^*))\Big)^{\widehat{M}}\rho_k(X^*)^{M-1}\\
        \geq&b_\mathbf{c}\Big(\min_{1\leq j\leq \widehat{M}}{\rm dist}(P_k(\mathbf{x}_n^*),P_{k}(\widehat{\mathbf{x}}_j^*))\Big)^{\widehat{M}}\rho_*^{M-1}.
    \end{align*}
Meanwhile, by the Cauchy-Schwarz inequality, we have
    \begin{align*}
        &\left|\mathcal{R}(\Psi_n;f-\widehat{f},g-\widehat{g})-\int_\Omega\Delta v_{\widehat{\theta}^*}\Psi_n{\rm d}\mathbf{x}\right|\\
        \leq&\|\Psi_n\|_{L^2(\partial\Omega)}\|g-\widehat{g}\|_{L^2(\partial\Omega)}+\|\partial_\nu\Psi_n\|_{L^2(\partial\Omega)}\|f-\widehat{f}\|_{L^2(\partial\Omega)}+\|\Delta v_{\widehat{\theta}^*}\|_{L^2(\Omega)}\|\Psi_n\|_{L^2(\Omega)}.
\end{align*}
By the definition of $\beta$, we have $|\Psi_n(\mathbf{x})|\leq \beta^{M+\widehat{M}-1}$ on $\overline{\Omega}$, and then
\begin{align*}
\|\Psi_n\|_{L^2(\partial\Omega)}
\leq |\partial\Omega|^{\frac{1}{2}}\beta^{M+\widehat{M}-1}\quad\mbox{and}
\quad
\|\Psi_n\|_{L^2(\Omega)}\leq |\Omega|^{\frac{1}{2}}\beta^{M+\widehat{M}-1}.
\end{align*}
Likewise, $|\partial_{x_k}\Psi_n(\mathbf{x})|\leq (M+\widehat{M}-1)\beta^{M+\widehat{M}-2}$, $|\partial_{x_{k+1}}\Psi_n(\mathbf{x})|\leq (M+\widehat{M}-1)\beta^{M+\widehat{M}-2}$ on $\overline{\Omega}$, and hence,
$\|\partial_\nu\Psi_n\|_{L^2(\partial\Omega)}\leq \sqrt{2}(M+\widehat{M}-1)|\partial\Omega|^{\frac{1}{2}}\beta^{M+\widehat{M}-2}.$
Thus, we have
\begin{align*}
    &\Big|\mathcal{R}(\Psi_n;f-\widehat{f},g-\widehat{g})-\int_\Omega\Delta v_{\widehat{\theta}^*}\Psi_n{\rm d}\mathbf{x}\Big| \\  \leq& C(\Omega,\beta,M,\widehat{M})(\|g-\widehat{g}\|_{L^2(\partial\Omega)}+\|f-\widehat{f}\|_{L^2(\partial\Omega)}+\|\Delta v_{\widehat{\theta}^*}\|_{L^2(\Omega)}).
\end{align*}
Then, taking the maximum over $ n$ and using the definition of $L(v_{\widehat{\theta}^*},\widehat{\mathbf{c}}^*,\widehat{X}^*)$ lead to
\begin{equation}
    \begin{aligned}
        &\max_{1\leq n\leq M}\min_{1\leq j\leq \widehat{M}}{\rm dist}\left(P_k(\mathbf{x}_n^*),P_{k}(\widehat{\mathbf{x}}_j^*)\right) \\
        \leq&\left[{C}(\Omega,\beta,M,\widehat{M})b_\mathbf{c}^{-1}\rho_*^{1-M}\big(\|g-\widehat{g}\|_{L^2(\partial\Omega)}+\|f-\widehat{f}\|_{L^2(\partial\Omega)}+\|\Delta v_{\widehat{\theta}^*}\|_{L^2(\Omega)}\big)\right]^{\tfrac{1}{\widehat{M}}}\\
        \leq& C(\Omega,\beta,M,\widehat{M},\boldsymbol{\sigma},b_\mathbf{c},\rho_*)L(v_{\widehat{\theta}^*},\widehat{\mathbf{c}}^*,\widehat{X}^*)^{\frac{1}{2\widehat{M}}}.
    \end{aligned}
    \end{equation}
Similarly, with $\widetilde{\Psi}_n:=\prod_{j=1}^{M}\left(P_k(\mathbf{x})-P_{k}^*(\mathbf{x}_j)\right)\prod_{j=1,j\neq n}^{\widehat{M}}(P_k(\mathbf{x})-P_{k}({\widehat{\mathbf{x}}_j^*}))$, we have
    \begin{equation}
    \begin{aligned}
        &\max_{1\leq n\leq \widehat{M}}\min_{1\leq j\leq M}{\rm dist}\left(P_{k}(\mathbf{x}_j^*),{P_{k}(\widehat{\mathbf{x}}_n^*)}\right) \\
        \leq&\left[{C}(\Omega,\beta,M,\widehat{M})b_\mathbf{c}^{-1}\rho_*^{1-\widehat{M}}\big(\|g-\widehat{g}\|_{L^2(\partial\Omega)}+\|f-\widehat{f}\|_{L^2(\partial\Omega)}+\|\Delta v_{\widehat{\theta}^*}\|_{L^2(\Omega)}\big)\right]^\frac{1}{M}\\
        \leq&C(\Omega,\beta,M,\widehat{M},\boldsymbol{\sigma},b_\mathbf{c},\rho_*)L(v_{\widehat{\theta}^*},\widehat{\mathbf{c}}^*,\widehat{X}^*)^{\frac{1}{2M}}.
    \end{aligned}
\end{equation}
Now using the definition of Hausdorff distance gives
\begin{equation}\label{eqn:project-Hausd}
d_H(\mathbf{P}_k(\widehat{X}^*), \mathbf{P}_k(X^*))\leq C(\Omega,\beta,M,\widehat{M},\boldsymbol{\sigma},b_\mathbf{c},\rho_*)L(v_{\widehat{\theta}^*},\widehat{\mathbf{c}}^*,\widehat{X}^*)^{\frac{1}{2\max\{M,\widehat{M}\}}},
\end{equation}
Further, direct computation yields
    $$\begin{aligned}
        \max_{1\leq n\leq M}\min_{1\leq j\leq \widehat{M}}{\rm dist}(\mathbf{x}_{n}^*,\widehat{\mathbf{x}}_{j}^*)\leq \sum_{k=1}^{d-1}\max_{1\leq n\leq M}\min_{1\leq j\leq \widehat{M}}{\rm dist}\left(P_{k}(\mathbf{x}_n^*),P_{k}(\widehat{\mathbf{x}}_j^*)\right),\\
        \max_{1\leq n\leq \widehat{M}}\min_{1\leq j\leq M}{\rm dist}(\mathbf{x}_{j}^*,\widehat{\mathbf{x}}_{n}^*)\leq \sum_{k=1}^{d-1}\max_{1\leq n\leq \widehat{M}}\min_{1\leq j\leq M}{\rm dist}\left(P_{k}(\mathbf{x}_j^*),P_{k}(\widehat{\mathbf{x}}_n^*)\right).
    \end{aligned}$$
That is, $d_H(\widehat{X}^*, X^*)\leq{\sum_{k=1}^{d-1}}d_H(\mathbf{P}_k(\widehat{X}^*), \mathbf{P}_k(X^*))$. The assertion follows from \eqref{eqn:project-Hausd}.
\end{proof}

By Theorem \ref{thm:disbound} and the Hall-Rado theorem in graph theory, we have the following error bound when $\widehat{M}=M$. The condition $\widehat M = M$ will be assume below.
\begin{theorem}\label{thm:disbound-1}
Let Assumption \ref{assump:separation} be fulfilled.
For two point configurations $\widehat{X}^*$ and $X^*$, if $\widehat{M}=M$, there exists a permutation $\pi$ of the set $\{1,2,\cdots,M\}$ such that
    \begin{equation}
        \max_{1\leq j\leq M}{\rm dist}(\widehat{\mathbf{x}}_{j}^*,\mathbf{x}_{\pi(j)}^*)\leq C(d-1)L(v_{\widehat{\theta}^*},\widehat{\mathbf{c}}^*,\widehat{X}^*)^{\frac{1}{2M}}.
    \end{equation}
\end{theorem}

\begin{theorem}[Hall-Rado \cite{Rado:1949}]\label{thm:hall-rado}
    Consider an even graph having $2M$ nodes $a_1, \ldots, a_M$ and $b_1, \ldots, b_M$, and connect some pairs $\left(a_i, b_j\right)$ such that for every $k \in\{1, \ldots, M\}$ and every subsequence $\left(a_{j_1}, \ldots, a_{j_k}\right)$ of $\left(a_1, \ldots, a_M\right)$, at least $k$ elements of the sequence $\left(b_1, \ldots, b_M\right)$ are connected to one of them. Then there exists a permutation $\pi$ of the integers $1, \ldots, M$ such that $a_j$ is connected to $b_{\pi(j)}$ for every $j$.
\end{theorem}

\begin{proof}[Proof of Theorem \ref{thm:disbound-1}]
The result is direct from Theorem \ref{thm:disbound}; see the proof of \cite[Theorem 3.2]{Badia2013}. We recap the proof for completeness. We connect the even graph $(\widehat{X}^*,X^*)$ in the way that $\widehat{\mathbf{x}}_i^*$ is connected to $\mathbf{x}_j^*$ if ${\rm dist}(\widehat{\mathbf{x}}_i^*, \mathbf{x}_j^*)\leq d_H(\widehat{X}^*,X^*)$ for $i,j=1,2,\cdots,M$. We claim that for every $k=1,2,\cdots,M$ and every subsequence $(\mathbf{x}_{j_1}^*,\mathbf{x}_{j_2}^*,\cdots,\mathbf{x}_{j_k}^*)$ of $X^*$, at least $k$ elements of $\widehat{\mathbf{x}}_i^*$ $(i=1,2,\cdots,M)$ are connected to one of them. Otherwise, there exists one subsequence $X_0^*=(\mathbf{x}_{j_1}^*,\mathbf{x}_{j_2}^*,\cdots,\mathbf{x}_{j_{k_0}}^*)$ of $X^*$ and at least one $\mathbf{x}_{j_t}^*$ in $X_0^*$, there holds ${\rm dist}(\mathbf{x}_{j_t}^*, \widehat{\mathbf{x}}_i^*)>d_H(\widehat{X}^*,X^*)$ for all $i=1,2,\cdots,M$, which leads to a contradiction with the definition of Hausdorff distance. Then the desired result follows directly from the Hall–Rado theorem.
\end{proof}

Next we bound the error of the recovered densities $\widehat{\mathbf{c}}^*$. \begin{theorem}\label{thm:densbound}
For the two configurations $(\mathbf{c}^*, X^*)$ and $(\widehat{\mathbf{c}}^*, \widehat{X}^*)$, let $\pi$ be the permutation given in Theorem \ref{thm:disbound-1}. Then under Assumption \ref{assump:separation}, there holds
\begin{equation}
        \max_{1\leq j\leq M}\big|\widehat{c}_j^*-c_{\pi(j)}^*\big|\leq C(\Omega,\beta,M,\boldsymbol{\sigma},b_\mathbf{c},B_\mathbf{c},\rho_*)L(v_{\widehat{\theta}^*},\widehat{\mathbf{c}}^*,\widehat{X}^*)^{\frac{1}{2M}}.
    \end{equation}
\end{theorem}
\begin{proof}
Fix $\mathbf{x}_0\in\Omega_\gamma$, and let $P_{1}(\mathbf{x}_0)$ be the projection of $\mathbf{x}_0$ onto the $x_1Ox_2$-complex plane.
For $n=0,1,\cdots,M-1$, define the harmonic functions
\begin{equation*}
\psi_n(\mathbf{x})=\beta^{-n}\left({P_{1}(\mathbf{x})-P_{1}(\mathbf{x}_0)}\right)^n=\beta^{-n}\left({x_1+{\rm i}x_2-P_{1}(\mathbf{x}_0)}\right)^n:= \overline{P}_1(\mathbf{x})^n.
\end{equation*}
Testing \eqref{equation-num-true} with $\psi_n$ and using  $\mathcal{R}(\psi_n)$ in \eqref{eq:R-formula}, with $\widehat f$ and $\widehat g$ defined in \eqref{eqn:hatfg}, yield
\begin{align}\label{eqn:R-density}
    \mathcal{R}(\psi_n;f,g)&=\sum_{j=1}^Mc_{\pi(j)}^*\psi_n(\mathbf{x}_{\pi(j)}^*)=\sum_{j=1}^Mc_{\pi(j)}^*\overline{P}_1(\mathbf{x}_{\pi(j)}^*)^n,\\
\label{eqn:Rhat-density}
    \mathcal{R}(\psi_n;\widehat{f},\widehat{g})&=\sum_{j=1}^{M}\widehat{c}_j^*\psi_n(\widehat{\mathbf{x}}_{j}^*)-\int_\Omega\Delta v_{\widehat{\theta}^*}\psi_n{\rm d}\mathbf{x}=\sum_{j=1}^M\widehat{c}_j^*\overline{P}_1(\widehat{\mathbf{x}}_{j}^*)^n-\int_\Omega\Delta v_{\widehat{\theta}^*}\psi_n{\rm d}\mathbf{x}.
\end{align}
Next we define two Vandermonde matrices
\begin{align*}
    V(P^*)&=\begin{pmatrix}
        1&1&\cdots&1\\
        \overline{P}_1(\mathbf{x}_{\pi(1)}^*)&\overline{P}_1(\mathbf{x}_{\pi(2)}^*)&\cdots&\overline{P}_1(\mathbf{x}_{\pi(M)}^*)\\
        \vdots&\vdots&\ddots&\vdots\\
        \overline{P}_1(\mathbf{x}_{\pi(1)}^*)^{M-1}&\overline{P}_1(\mathbf{x}_{\pi(2)}^*)^{M-1}&\cdots&\overline{P}_1(\mathbf{x}_{\pi(M)}^*)^{M-1}
    \end{pmatrix},\\    V(\widehat{P}^*)&=\begin{pmatrix}
        1&1&\cdots&1\\
        \overline{P}_1(\widehat{\mathbf{x}}_{1}^*)&\overline{P}_1(\widehat{\mathbf{x}}_{2}^*)&\cdots&\overline{P}_1(\widehat{\mathbf{x}}_{M}^*)\\
        \vdots&\vdots&\ddots&\vdots\\
        \overline{P}_1(\widehat{\mathbf{x}}_{1}^*) ^{M-1}&\overline{P}_1(\widehat{\mathbf{x}}_{2}^*)^{M-1}&\cdots&\overline{P}_1(\widehat{\mathbf{x}}_{M}^*)^{M-1}
    \end{pmatrix},
\end{align*}
and four vectors in $\mathbb{R}^M$, i.e.,
\begin{equation*}
\begin{split}
\mathcal{R}^*&=(\mathcal{R}(\psi_0;f,g),\mathcal{R}(\psi_1;f,g),\cdots,\mathcal{R}(\psi_{M-1};f,g))^{\rm T},\\
 \widehat{\mathcal{R}}&=(\mathcal{R}(\psi_0;\widehat{f},\widehat{g}),\mathcal{R}(\psi_1;\widehat{f},\widehat{g}),\cdots,\mathcal{R}(\psi_{M-1};\widehat{f},\widehat{g}))^{\rm T},\\
 \widetilde{\mathbf{c}}^*&=(c_{\pi(1)}^*, c_{\pi(2)}^*, \cdots, c_{\pi(M)}^*)^{\rm T}, ~~ \text{and}\\
 \mathcal{Q} & =  (\int_\Omega\Delta v_{\widehat{\theta}^*}\psi_0{\rm d}\mathbf{x},\int_\Omega\Delta v_{\widehat{\theta}^*}\psi_1{\rm d}\mathbf{x},\cdots,\int_\Omega\Delta v_{\widehat{\theta}^*}\psi_{M-1}{\rm d}\mathbf{x} )^{\rm T}
\end{split}
\end{equation*}

Then we can rewrite \eqref{eqn:R-density} and \eqref{eqn:Rhat-density} into
$V(P^*)\widetilde{\mathbf{c}}^*=\mathcal{R}^*$  and $ V(\widehat{P}^*)\widehat{\mathbf{c}}^*=\widehat{\mathcal{R}}-\mathcal{Q}$.
Taking their difference yields
\begin{align*}
\widetilde{\mathbf{c}}^*-\widehat{\mathbf{c}}^*=&V(P^*)^{-1}\mathcal{R}^*-\widehat{\mathbf{c}}^*=V(P^*)^{-1}(\mathcal{R}^*-\widehat{\mathcal{R}})+V(P^*)^{-1}\widehat{\mathcal{R}}-V(P^*)^{-1}V(P^*)\widehat{\mathbf{c}}^*\\
=& V(P^*)^{-1}\big[( \mathcal{R}^*-\widehat{\mathcal{R}})+(V(\widehat{P}^*)\widehat{\mathbf{c}}^*-V(P^*)\widehat{\mathbf{c}}^*)+\mathcal{Q}\big],
\end{align*}
and consequently,
\begin{align}
\|\widetilde{\mathbf{c}}^*-\widehat{\mathbf{c}}^*\|_\infty
\leq\|V(P^*)^{-1}\|_\infty\big(\|\mathcal{R}^*-\widehat{\mathcal{R}}\|_\infty+\|V(P^*)\widehat{\mathbf{c}}^*-V(\widehat{P}^*)\widehat{\mathbf{c}}^*\|_\infty+\|\mathcal{Q}\|_\infty\big).\label{ineq:densbound}
\end{align}
Meanwhile, \cite[Theorem 1]{Gautschi1962OnIO} gives
\begin{equation}\label{eqn:densbound-v}
    \|V(P^*)^{-1}\|_\infty\leq\max_{1\leq j\leq M}\prod_{i=1,i\neq j}^M\dfrac{1+\left|\overline{P}_1(\mathbf{x}_{\pi(i)}^*)\right|}{\left|\overline{P}_1(\mathbf{x}_{\pi(j)}^*)-\overline{P}_1(\mathbf{x}_{\pi(i)}^*)\right|}\leq\left(\dfrac{2\beta}{\rho_*}\right)^{M-1}.
\end{equation}
Since $\mathbf{x}_0\in \Omega_\gamma$, we have $|\psi_n(\mathbf{x})|\leq1$, and $|\partial_{x_1}\psi_n(\mathbf{x})|\leq n\beta^{-1}, |\partial_{x_2}\psi_n(\mathbf{x})|\leq n\beta^{-1}$, and $|\partial_{\nu}\psi_n(\mathbf{x})|\leq\sqrt{2}n\beta^{-1}$. This leads to
\begin{align}
        \|\mathcal{R}^*-\widehat{\mathcal{R}}\|_\infty&\leq\max_{0\leq n\leq M-1}\!\!\!\|\psi_n\|_{L^2(\partial\Omega)}\|g-\widehat{g}\|_{L^2(\partial\Omega)}+\max_{0\leq n\leq M-1}\!\!\|\partial_\nu\psi_n\|_{L^2(\partial\Omega)}\|f-\widehat{f}\|_{L^2(\partial\Omega)}\nonumber\\
        &\leq |\partial\Omega|^{\frac{1}{2}}\|g-\widehat{g}\|_{L^2(\partial\Omega)}+\sqrt{2}(M-1)\beta^{-1}|\partial\Omega|^{\frac{1}{2}}\|f-\widehat{f}\|_{L^2(\partial\Omega)}.\label{eqn:densbound-1}
    \end{align}
Likewise, we have
\begin{align*}        &\|V(P^*)\widehat{\mathbf{c}}^*-V(\widehat{P}^*)\widehat{\mathbf{c}}^*\|_\infty=\Big\|\sum_{j=1}^M \widehat{c}_j^*\left(\psi_n(\mathbf{x}_{\pi(j)}^*)-\psi_n(\widehat{\mathbf{x}}_j^*)\right)\Big\|_{\infty}\\
        \leq&\sqrt{2}\beta^{-1}\sum_{j=1}^{M}|\widehat{c}_j^*|\max_{0\leq n\leq M-1}\|\nabla\psi_n\|_{L^\infty(\Omega;\mathbb{C}^2)}\max_{1\leq j\leq M}\|P_{1}(\mathbf{x}_{\pi(j)}^*)-P_{1}(\widehat{\mathbf{x}}_j^*)\|\\
        \leq& \sqrt{2}\beta^{-2}M(M-1)B_{\mathbf{c}}\max_{1\leq j\leq M}\|\mathbf{x}_{\pi(j)}^*-\widehat{\mathbf{x}}_j^*\|,\\    &\|\mathcal{Q}\|_\infty\leq\|\Delta v_{\widehat{\theta}^*}\|_{L^2(\Omega)}\max_{0\leq n\leq M-1}\|\psi_n\|_{L^2(\Omega)}\leq\|\Delta v_{\widehat{\theta}^*}\|_{L^2(\Omega)}.
\end{align*}
Combining the preceding estimates yields with $C=C(\Omega,\beta,M,B_\mathbf{c},\rho_*,\boldsymbol{\sigma},d)$,
\begin{align*}
    \|\widetilde{\mathbf{c}}^*-\widehat{\mathbf{c}}^*\|_\infty
    \leq&C\big(\|g-\widehat{g}\|_{L^2(\partial\Omega)}+\|f-\widehat{f}\|_{L^2(\partial\Omega)}+\|\Delta v_{\widehat{\theta}^*}\|_{L^2(\Omega)}+\max_{1\leq j\leq M}\|\mathbf{x}_{\pi(j)}^*-\widehat{\mathbf{x}}_j^*\|\big)\\
    \leq&C\big(L(v_{\widehat{\theta}^*},\widehat{\mathbf{c}}^*,\widehat{X}^*)^{\frac{1}{2}}+\max_{1\leq j\leq M}\|\mathbf{x}_{\pi(j)}^*-\widehat{\mathbf{x}}_j^*\|\big).
\end{align*}
 This estimate and Theorem \ref{thm:disbound-1} imply the desired result.
\end{proof}

Next we analyze the error of the regular part $v_{\widehat\theta^*}$ in terms of the population loss.
\begin{theorem}\label{thm:errorv}
Let $v_{\widehat{\theta}^*}$ be defined in \eqref{eqn:min-emp}, and $v^*$ the solution to \eqref{problem-v}. Then under Assumption \ref{assump:separation},
\begin{equation}
    \|v^*-v_{\widehat{\theta}^*}\|_{L^2(\Omega)}\leq C(\Omega,\beta,M,\boldsymbol{\sigma},b_\mathbf{c},B_\mathbf{c},\rho_*,\gamma,d)L(v_{\widehat{\theta}^*},\widehat{\mathbf{c}}^*,\widehat{X}^*)^{\frac{1}{2M}}.
\end{equation}
\end{theorem}
\begin{proof}
Let $\pi$ be the permutation given in Theorem \ref{thm:disbound-1}, let $e=v^*-v_{\widehat{\theta}^*}$, and let
\begin{equation}\label{def:zeta}
    \zeta=\sum_{j=1}^M(c_{\pi(j)}^*\Phi_{\mathbf{x}_{\pi(j)}^*}-\widehat{c}_j^*\Phi_{\widehat{\mathbf{x}}_j^*})=\sum_{j=1}^M\big(c_{\pi(j)}^*-\widehat{c}_j^*\big)\Phi_{\mathbf{x}_{\pi(j)}^*}+\sum_{j=1}^M\widehat{c}_j^*\big(\Phi_{\mathbf{x}_{\pi(j)}^*}-\Phi_{\widehat{\mathbf{x}}_j^*}\big).
\end{equation}
Then with $\widehat f$ and $\widehat g$ in \eqref{eqn:hatfg}, we have
\begin{align}
L(v_{\widehat{\theta}^*},\widehat{\mathbf{c}}^*,\widehat{X}^*)
=&\|\Delta v_{\widehat{\theta}^*}\|_{L^2(\Omega)}^2+\sigma_d\|\widehat f-f\|_{L^2(\partial\Omega)}^2+\sigma_n\|\widehat g-g\|_{L^2(\partial\Omega)}^2\nonumber\\
    =&\|\Delta e\|_{L^2(\Omega)}^2+\sigma_d\|e+\zeta\|_{L^2(\partial\Omega)}^2+\sigma_n\|\partial_\nu e+\partial_\nu\zeta\|_{L^2(\partial\Omega)}^2.\label{eqn:L-vhatstar}
\end{align}
Now we define the harmonic extension $\eta$ of $e$ by
\begin{equation*}
\left\{\begin{aligned}-\Delta \eta &=0 , &&\text {in } \Omega,\\
\eta &=e , &&\text {on } \partial\Omega.
\end{aligned}\right.
\end{equation*}
The following elliptic regularity estimate holds \cite[Theorem 4.2, p. 870]{Berggren:2004}
\begin{equation}\label{eqn:stab-hm}
\|  \eta  \|_{L^2(\Omega)} \le C(\Omega)  \| e \|_{L^2(\partial\Omega)} .
\end{equation}
Let $\tilde e = e - \eta$. Then it satisfies
\begin{equation*}
\left\{\begin{aligned}
-\Delta \tilde e &=  - \Delta e  , &&\text {in } \Omega,\\
\tilde e &= 0 , && \text {on } \partial\Omega.
\end{aligned}\right.
\end{equation*}
Since $\Delta e \in L^2(\Omega)$, the Poincar\'e inequality and the standard energy argument imply
$$ \|\tilde e\|_{L^2(\Omega)} \le C(\Omega) \| \nabla \tilde e  \|_{L^2(\Omega)} \le C(\Omega) \|\Delta e\|_{L^2(\Omega)}.$$
This, the triangle inequality, and the stability estimate \eqref{eqn:stab-hm} lead to
\begin{align*}
    \|e\|_{L(\Omega)}^2\leq& 2(\|\tilde e\|_{L^2(\Omega)}^2+\|\eta\|_{L^2(\Omega)}^2)
    \leq C(\Omega,\boldsymbol{\sigma})\big(\|\Delta e\|_{L^2(\Omega)}^2+\sigma_d\|e\|_{L^2(\partial\Omega)}^2\big)\\
    \leq&C(\Omega,\boldsymbol{\sigma})(\|\Delta e\|_{L^2(\Omega)}^2+\sigma_d\|e\|_{L^2(\partial\Omega)}^2+\sigma_n\|\partial_\nu e\|_{L^2(\partial\Omega)}^2).
\end{align*}
This estimate, the triangle inequality and the identity \eqref{eqn:L-vhatstar} imply
\begin{align}
    \|e\|_{L(\Omega)}^2
\leq&C(\Omega,\boldsymbol{\sigma})\left(L(v_{\widehat{\theta}^*},\widehat{\mathbf{c}}^*,\widehat{X}^*)+\sigma_d\|\zeta\|_{L^2(\partial\Omega)}^2+\sigma_n\|\partial_\nu\zeta\|_{L^2(\partial\Omega)}^2\right).\label{ineq:eL2}
\end{align}
Meanwhile, by the definition of $\Phi(\mathbf{x})$ in \eqref{fundamental-s},
\begin{align*}\max_{1\leq j\leq M}\left\|\Phi_{\mathbf{x}_{j}^*}\right\|_{L^2(\partial\Omega)}&\leq
c_d|\partial\Omega|^{\frac{1}{2}}\left\{
\begin{aligned}
	\max\{|\log\gamma|,|\log\beta|\},&\quad d=2,\\
		\gamma^{2-d},&\quad d\geq3,\\
	\end{aligned}\right.\\
\max_{1\leq j\leq M}\|\partial_\nu\Phi_{\mathbf{x}_{j}^*}\|_{L^2(\partial\Omega)}&\leq c_d'|\partial\Omega|^{\frac{1}{2}}
	\gamma^{1-d},\\
\max_{1\leq j\leq M}\|\Phi_{\mathbf{x}_{\pi(j)}^*}-\Phi_{\widehat{\mathbf{x}}_j^*}\|_{L^2(\partial\Omega)}&\leq c_d'|\partial\Omega|^{\frac{1}{2}}\gamma^{1-d}\|\mathbf{x}_{\pi(j)}^*-\widehat{\mathbf{x}}_j^*\|,\\
\max_{1\leq j\leq M}\left\|\partial_\nu\Phi_{\mathbf{x}_{\pi(j)}^*}-\partial_\nu\Phi_{\widehat{\mathbf{x}}_j^*}\right\|_{L^2(\partial\Omega)}&\leq c_d'|\partial\Omega|^{\frac{1}{2}}(d-1)\gamma^{-d}\|\mathbf{x}_{\pi(j)}^*-\widehat{\mathbf{x}}_j^*\|,
\end{align*}
with $c_d'={(d\alpha(d))}^{-1}$.
Next we bound $\zeta$ and $\partial_\nu\zeta$. From the splitting \eqref{def:zeta}
and by combining the preceding estimates with Theorems \ref{thm:disbound-1} and \ref{thm:densbound}, we obtain
\begin{align*}
  \|\zeta\|_{L^2(\Omega)}  \leq&M\max_{1\leq j\leq M}|\widehat{c}_j^*-c_{\pi(j)}^*|\max_{1\leq j\leq M}\|\Phi_{\mathbf{x}_{j}^*}\|_{L^2(\partial\Omega)}\\
    &+M\max_{1\leq j\leq M}|\widehat{c}_j^*|\max_{1\leq j\leq M}\|\Phi_{\mathbf{x}_{\pi(j)}^*}-\Phi_{\widehat{\mathbf{x}}_j^*}\|_{L^2(\partial\Omega)}\leq CL(v_{\widehat{\theta}^*},\widehat{\mathbf{c}}^*,\widehat{X}^*)^{\frac{1}{2M}},\\
    \|\partial_\nu\zeta\|_{L^2(\partial\Omega)}
    \leq&M\max_{1\leq j\leq M}|\widehat{c}_j^*-c_{\pi(j)}^*|\max_{1\leq j\leq M}\|\partial_\nu\Phi_{\mathbf{x}_{j}^*}\|_{L^2(\partial\Omega)}\\
     &+M\max_{1\leq j\leq M}\left|\widehat{c}_j^*\right|\max_{1\leq j\leq M}\|\partial_\nu\Phi_{\mathbf{x}_{\pi(j)}^*}-\partial_\nu\Phi_{\widehat{\mathbf{x}}_j^*}\|_{L^2(\partial\Omega)}\leq CL(v_{\widehat{\theta}^*},\widehat{\mathbf{c}}^*,\widehat{X}^*)^{\frac{1}{2M}}.
\end{align*}
Then substituting these two estimates into \eqref{ineq:eL2} yields the desired result.
\end{proof}

Finally, we gives an error bound on the approximation of the solution $u^*$.
\begin{corollary}\label{cor:erroru}
Let $u_{\widehat{\theta}^*}=v_{\widehat{\theta}^*} + \sum_{j=1}^M\widehat{c}_j^*\Phi_{\widehat{\mathbf{x}}_j^*}$, and $u^*$ be the solution of \eqref{problem}. Then for any $p<\frac{d}{d-1}$, there holds with $C=C(\Omega,\beta,M,\boldsymbol{\sigma}, b_\mathbf{c},B_\mathbf{c},\rho_*, \gamma,p,d)$,
\begin{equation*}
\|u^*-u_{\widehat{\theta}^*}\|_{L^p(\Omega)} \leq CL(v_{\widehat{\theta}^*},\widehat{\mathbf{c}}^*,\widehat{X}^*)^{\frac{1}{2M}}.
\end{equation*}
\end{corollary}
\begin{proof} With $\zeta$ defined in \eqref{def:zeta}, by H{\"o}lder's inequality, we have
    $$\|u^*-u_{\widehat{\theta}^*}\|_{L^p(\Omega)}\leq \|v^*-v_{\widehat{\theta}^*}\|_{L^p(\Omega)}+\|\zeta\|_{L^p(\Omega)}\leq |\Omega|^{\frac{2-p}{2p}}\|v^*-v_{\widehat{\theta}^*}\|_{L^2(\Omega)}+\|\zeta\|_{L^p(\Omega)}.$$
    So it suffices to bound $\|\zeta\|_{L^{p}(\Omega)}$. Since $\Phi_{\mathbf{x}_j^*}\in W^{1,p}(\Omega)$ \cite[Theorem 1.1, p. 3]{GruterWidman:1982} and $F(\mathbf{y})=\|\nabla\Phi_{\mathbf{y}}(\mathbf{x})\|_{L^{p}(\Omega)}=(\int_\Omega(c_d'|\mathbf{x}-\mathbf{y}|^{1-d})^p{\rm d}\mathbf{x})^{\frac{1}{p}}$ attains its maximum in $\Omega_\gamma$,
\begin{align*}    &\|\zeta\|_{L^p(\Omega)}\leq\Big\|\sum_{j=1}^M\left(c_{\pi(j)}^*-\widehat{c}_j^*\right)\Phi_{\mathbf{x}_{\pi(j)}^*}\Big\|_{L^p(\Omega)}+\Big\|\sum_{j=1}^M\widehat{c}_j^*\left(\Phi_{\mathbf{x}_{\pi(j)}^*}- \Phi_{\widehat{\mathbf{x}}_j^*}\right)\Big\|_{L^p(\Omega)}\\
    \leq&M\Big(\max_{1\leq j\leq M}|\widehat{c}_j^*-c_{\pi(j)}^*|\max_{1\leq j\leq M}\|\Phi_{\mathbf{x}_{j}^*}\|_{L^p(\Omega)}+\max_{1\leq j\leq M}|\widehat{c}_j^*|\max_{1\leq j\leq M}\|\Phi_{\mathbf{x}_{\pi(j)}^*}-\Phi_{\widehat{\mathbf{x}}_j^*}\|_{L^p(\Omega)}\Big)\\
    \leq&C\Big(\max_{1\leq j\leq M}|\widehat{c}_j^*-c_{\pi(j)}^*|+\max_{1\leq j\leq M}\|\widehat{\mathbf{x}}_{j}^*-\mathbf{x}_{\pi(j)}^*\|\Big)\leq CL(v_{\widehat{\theta}^*},\widehat{\mathbf{c}}^*,\widehat{X}^*)^{\frac{1}{2M}}.
\end{align*}
This completes the proof of the corollary.
\end{proof}

\subsection{Generalization error analysis of the empirical loss}
The analysis in Section \ref{ssec:stability} focuses on the population loss, evaluated at the minimizer of the empirical loss $\widehat {L}^\delta(\theta,\mathbf{c},X)$. In practice, we only minimize the empirical loss $\widehat {L}^\delta(\theta,\mathbf{c},X)$, and do not have access to the population loss  $L(\theta,\mathbf{c},X)$. The analysis of the error between the empirical and population losses is commonly known as generalization error analysis in statistical learning theory \cite{AnthonyBartlett:1999}, which is the main goal of this part. First we establish an important error splitting, which involves three components: the approximation error $\inf_{v_\theta\in\mathcal{A}}\|v_\theta-v^*\|^2_{H^2(\Omega)} $, the statistical error $ \sup_{(v_\theta,\mathbf{c},X)\in\mathbb{A}}|\widehat{L}^\delta(v_\theta,\mathbf{c},X)-L^\delta(v_\theta,\mathbf{c},X)|$ and the noise level $\delta^2$. The first two components will be  analyzed below.
\begin{lemma}\label{lemma:errordecom}
Let $(\widehat{\theta}^*,\widehat{\mathbf{c}}^*,\widehat{X}^*)$ minimize the loss \eqref{eqn:loss-emp} over the set $\mathbb{A}$. Then we have
\begin{align*}
L(v_{\widehat{\theta}^*},\widehat{\mathbf{c}}^*,\widehat{X}^*)
\leq
4\sup_{(v_\theta,\mathbf{c},X)\in\mathbb{A}}\left|\widehat{L}^\delta(v_\theta,\mathbf{c},X)-L^\delta(v_\theta,\mathbf{c},X)\right|+C(\boldsymbol{\sigma})\Big(\inf_{v_\theta\in\mathcal{A}}\|v_\theta-v^*\|^2_{H^2(\Omega)}+\delta^2\Big).\nonumber
\end{align*}
\end{lemma}
\begin{proof}
Let $F(\mathbf{c},X)=\sum_{j=1}^Mc_j\Phi_{\mathbf{x}_j}$.
By the triangle inequality, for any $(v,\mathbf{c},X)\in H^2(\Omega)\times\mathbb{R}^M\times\mathbb{R}^{d\times M}$, we have
\begin{align*}
&L^\delta(v,\mathbf{c},X)=\|\Delta v\|_{L^2(\Omega)}^2+\sigma_d\|v-f^\delta+F\|_{L^2(\partial\Omega)} ^2+\sigma_n\|\partial_\nu v-g^\delta+\partial_\nu F\|_{L^2(\partial\Omega)}^2\\
\leq& \|\Delta v\|_{L^2(\Omega)}^2+2\sigma_d\|v-f+F\|_{L^2(\partial\Omega)}^2+2\sigma_n\|\partial_\nu v-g+\partial_\nu F\|_{L^2(\partial\Omega)}^2+2(\sigma_d+\sigma_n)\delta^2\\
\leq&2L(v,\mathbf{c},X)+2(\sigma_d+\sigma_n)\delta^2.
\end{align*}
Similarly, we have $L(v,\mathbf{c},X)\leq2L^\delta(v,\mathbf{c},X)+2(\sigma_d+\sigma_n)\delta^2.$
These inequalities and the minimizing property of $(\widehat{\theta}^*,\widehat{\mathbf{c}}^*,\widehat{X}^*)$ to the loss $\widehat{L}^\delta$ imply that for any $v_\theta\in\mathcal{A}$,
    \begin{align} &L(v_{\widehat{\theta}^*},\widehat{\mathbf{c}}^*,\widehat{X}^*)
        \leq2L^\delta (v_{\widehat{\theta}^*},\widehat{\mathbf{c}}^*,\widehat{X}^*)+2(\sigma_d+\sigma_n)\delta^2\nonumber\\
        =& 2L^\delta (v_{\widehat{\theta}^*},\widehat{\mathbf{c}}^*,\widehat{X}^*)
        -2\widehat{L}^\delta (v_{\widehat{\theta}^*},\widehat{\mathbf{c}}^*,\widehat{X}^*)
        +2\widehat{L}^\delta (v_{\widehat{\theta}^*},\widehat{\mathbf{c}}^*,\widehat{X}^*)
        -2\widehat{L}^\delta (v_\theta,\mathbf{c}^*,X^*)\nonumber\\
        &+2\widehat{L}^\delta (v_\theta,\mathbf{c}^*,X^*)
        -2L^\delta (v_\theta,\mathbf{c}^*,X^*)+2L^\delta (v_\theta,\mathbf{c}^*,X^*)+2(\sigma_d+\sigma_n)\delta^2\nonumber\\
        \leq& 4\sup_{(v_\theta,\mathbf{c},X)\in\mathbb{A}}\left|\widehat{L}^\delta(v_\theta,\mathbf{c},X)-L^\delta(v_\theta,\mathbf{c},X)\right|+2L^\delta (v_\theta,\mathbf{c}^*,X^*)+2(\sigma_d+\sigma_n)\delta^2\nonumber\\
        \leq& 4\sup_{(v_\theta,\mathbf{c},X)\in\mathbb{A}}\left|\widehat{L}^\delta(v_\theta,\mathbf{c},X)-L^\delta(v_\theta,\mathbf{c},X)\right|+4L(v_\theta,\mathbf{c}^*,X^*)+6(\sigma_d+\sigma_n)\delta^2. \label{ineq:L-decom}
\end{align}
Using the facts $v^*|_{\partial\Omega}=f-F(\mathbf{c}^*,X^*)$ and $\partial_\nu v^*|_{\partial\Omega}=g-\partial_\nu F(\mathbf{c}^*,X^*)$, and the trace theorem, we deduce
    \begin{align*}
        &L(v_\theta,\mathbf{c}^*,X^*)=\|\Delta v_\theta\|_{L^2(\Omega)}^2 + \sigma_d\|v_\theta-v^*\|_{L^2(\partial\Omega)}^2+\sigma_n\|\partial_\nu v_\theta-\partial_\nu v^*\|_{L^2(\partial\Omega)}^2\nonumber\\
        =&\|\Delta v_\theta- \Delta v^*\|_{L^2(\Omega)}^2+\sigma_d\|v_\theta-v^*\|_{L^2(\partial\Omega)}^2+\sigma_n\|\partial_\nu v_\theta-\partial_\nu v^*\|_{L^2(\partial\Omega)}^2
        \leq C(\boldsymbol{\sigma})\|v_\theta-v^*\|_{H^2(\Omega)}^2.
    \end{align*}
Combining this with \eqref{ineq:L-decom} and taking the infimum over $v_\theta\in\mathcal{A}$ yield the result.
\end{proof}

Now we analyze the approximation and statistical errors separately. We make the following regularity assumption.  Condition (i) and the standard elliptic regularity theory ensure that problem \eqref{problem-v} admits a solution $v\in H^3(\Omega)$, which enables bounding the approximation error quantitatively. Condition (ii) is needed for bounding the Rademacher complexity of the function classes in Lemma \ref{lem:fcn-Lip}.
\begin{assumption}\label{assump:regularity} {\rm(i)} $f\in H^{\frac{5}{2}}(\partial\Omega)$ and $g\in H^{\frac{3}{2}}(\partial\Omega)$, and {\rm(ii)}  $f^\delta,g^\delta\in L^\infty(\partial\Omega)$.
\end{assumption}
 We employ the following approximation result \cite[Proposition 4.8]{GuhringRaslan:2021}. The notation $(s=2)$ equals to 1 if $s=2$ and zero otherwise.
\begin{lemma}\label{lem:tanh-approx}
Let $s\in\mathbb{N}\cup\{0\}$ and $p\in[1,\infty]$ be fixed, and $v\in W^{k,p}(\Omega)$ with $k\geq  s+1$. Then for any tolerance $\epsilon>0$, there exists at least one $v_\theta$ of depth $\mathcal{O}\big(\log(d+k)\big)$, with $|\theta|_{\ell^0}$ bounded by $\mathcal{O}\big(\epsilon^{-\frac{d}{k-s-\mu (s=2)}}\big)$ and  $|\theta|_{\ell^\infty}$ bounded by $\mathcal{O}(\epsilon^{-2-\frac{2(d/p+d+s+\mu(s=2))+d/p+d}{k-s-\mu (s=2)}})$, where $\mu>0$ is arbitrarily small, such that
$\|v-v_\theta\|_{W^{s,p}(\Omega)} \leq \epsilon$.
\end{lemma}

Under Assumption \ref{assump:regularity}, the regular part $v^*$ belongs to $H^3(\Omega)$. Thus, for any $\epsilon>0$, by Lemma \ref{lem:tanh-approx} with $k=3$ and $s=2$, there exists
\begin{equation}\label{eqn:NN-set}
v_\theta \in \mathcal{W}_\epsilon=:\mathcal{N}_{\varrho}(C, C\epsilon^{-\frac{d}{1-\mu}},C\epsilon^{-2-\frac{9d/2+4+2\mu}{1-\mu}}),
\end{equation}
such that $ \|w^*-v_\theta  \|_{H^2(\Omega)} \le \epsilon$.
Next, we bound the statistical error $$\mathcal{E}_{stat}=\sup_{(v_\theta,\mathbf{c},X)\in\mathcal{W}_\epsilon\times I_\mathbf{c}^M\times \Omega_\gamma^M}|\widehat{L}^\delta(v_\theta,\mathbf{c},X)-L^\delta(v_\theta,\mathbf{c},X)|,$$
which arises from approximating the integrals by Monte Carlo. Note that the discretization error is a random variable. In the statistical learning theory, such error is often bounded using Rademacher complexity. By the triangle inequality, we have
\begin{align}
\mathcal{E}_{stat}
 \le& \sup_{v_\theta\in \mathcal{W}_\epsilon} |\Omega|\Big| \frac{1}{N_r} \sum_{i=1}^{N_r} h_r(X_i;v_\theta) - \mathbb{E}_X (h_r(X;v_\theta)) \Big|\nonumber\\
&+\sup_{(v_\theta,\mathbf{c},X)\in\mathcal{W}_\epsilon\times I_\mathbf{c}^M\times \Omega_\gamma^M}  \sigma_d |\partial\Omega| \Big|  \frac{1}{N_d} \sum_{j=1}^{N_d} h_d(Y_j;v_\theta,\mathbf{c},X) - \mathbb{E}_Y (h_d(Y;v_\theta,\mathbf{c},X))  \Big|\nonumber\\
&+\sup_{(v_\theta,\mathbf{c},X)\in\mathcal{W}_\epsilon\times I_\mathbf{c}^M\times \Omega_\gamma^M}  \sigma_n |\partial\Omega| \Big|  \frac{1}{N_n} \sum_{k=1}^{N_n} h_n(Y_k;v_\theta,\mathbf{c},X) - \mathbb{E}_Y (h_n(Y;v_\theta,\mathbf{c},X))  \Big|,\label{eqn:err-sta}
 \end{align}
with
$h_r(\mathbf{x};v_\theta) = |\Delta v_\theta(\mathbf{x})|^2$ for $\mathbf{x}\in \Omega$, $h_d(\mathbf{\mathbf{y}};v_\theta,\mathbf{c},X) = (v_\theta-f^\delta+\sum_{j=1}^Mc_j\Phi_{\mathbf{x}_j})^2(\mathbf{y})$ and $h_n(\mathbf{y};v_\theta,\mathbf{c},X) = (\partial_\nu v_\theta-g^\delta+\sum_{j=1}^Mc_j\partial_\nu\Phi_{\mathbf{x}_j})^2(\mathbf{y}) $ for $\mathbf{y}\in \partial\Omega$. We define three function classes
$\mathcal{H}_r = \{h_r(v_\theta): v_\theta \in  \mathcal{W}_\epsilon\}$, $ \mathcal{H}_d  = \{h_d(v_\theta,\mathbf{c},X): v_\theta \in  \mathcal{W}_\epsilon, \mathbf{c}\in I_{\mathbf{c}}^M, X\in\Omega_\gamma^M \}$ and $\mathcal{H}_n  = \{h_n(v_\theta,\mathbf{c},X): v_\theta \in  \mathcal{W}_\epsilon, \mathbf{c}\in I_{\mathbf{c}}^M, X\in\Omega_\gamma^M\}$.

Using the technical estimates in Appendix \ref{app:est}, we have the following bound on $\mathcal{E}_{stat}$.
\begin{theorem}\label{thm:err-stat}
Under Assumption \ref{assump:regularity}, for small $\tau $, with probability at least $1- 3 \tau$,
$$\mathcal{E}_{stat}\leq C(e_r + \max\{\sigma_d,\sigma_n\} e_b),$$
with  $C={C}(d,\|f^\delta\|_{L^\infty(\Omega)},\|g^\delta\|_{L^\infty(\Omega)},B_\mathbf{c},\gamma,\beta)$, and the errors $e_r$ and $e_b$  defined by
\begin{align*}
  e_r & \le C \frac{ B_\theta^{4L} N_\theta^{4L-4} \big(N_\theta^\frac12\big( \log^\frac12 B_\theta + \log^\frac12 N_\theta + \log^\frac12 N_r)  + \log^\frac12 \frac{1}{\tau}\big)}{\sqrt{N_r}}, \\
 e_b & \le C \frac{B_\theta^{2L} N_\theta^{\max(2L-2,2)} \big(N_\theta^{\frac12} (\log^\frac12 B_\theta+ \log^\frac12 N_\theta + \log^\frac12 N_b) + \log^\frac12\frac{1}{\tau} \big)}{\sqrt{N_b}}.
 \end{align*}
\end{theorem}

\begin{proof}
Fix $m \in \mathbb{N }$, $B_\theta \in [1, \infty)$, $\epsilon \in  (0,1)$,
and $\mathbb{B}_{B_\theta} := \{x\in\mathbb{R}^m:\ |x|_{\ell^\infty}\leq B_\theta\}$. By \cite[Prop. 5]{CuckerSmale:2002}, $\mathcal{C}(\mathbb{B}_{B_\theta},|\cdot|_{\ell^\infty},\epsilon)$ (cf. Def. \ref{def:cover}) is bounded by
$ \log \mathcal{C}(\mathbb{B}_{B_\theta},|\cdot|_{\ell^\infty},\epsilon)\leq m\log (4B_\theta\epsilon^{-1}).$
Lemmas \ref{lem:NN-Lip} and \ref{lem:fcn-Lip} imply with $\Lambda_{r} = C N_\theta L^3 W^{5L-5}B_\theta^{5L-3}$,
\begin{align*}
 \log \mathcal{C}(\mathcal{H}_r,\|\cdot\|_{L^{\infty}(\Omega)},\epsilon)&\leq
\log \mathcal{C}(\Theta ,|\cdot|_{\ell^\infty},\Lambda_{r}^{-1}\epsilon)  \leq CN_\theta \log(4B_\theta \Lambda_r\epsilon^{-1}).
\end{align*}
By Lemma \ref{lem:fcn-Lip},
we have $M_{\mathcal{H}_r} = C L^2 W^{4L-4}B_\theta^{4L}$. Then by setting $s=n^{-\frac12}$ in Lemma \ref{lem:Dudley} and using the facts $1\leq B_\theta$, $1\leq L$ and $1\leq W \leq N_\theta$, $1\leq L\leq C\log (d+2)$ (for $k=3$), cf. Lemma \ref{lem:tanh-approx}, we can bound the Rademacher complexity $\mathfrak{R}_n(\mathcal{H}_r) $ by
\begin{align*}
\mathfrak{R}_n(\mathcal{H}_r)
 \leq&4n^{-\frac12} +12n^{-\frac12}\int^{M_{\mathcal{H}_r}}_{n^{-\frac12}}{\big(CN_\theta \log(4B_\theta \Lambda_{r}\epsilon^{-1})\big)}^{\frac12}\ {\rm d}\epsilon\\
 \leq& 4n^{-\frac12}+12n^{-\frac12}M_{\mathcal{H}_r}\big(CN_\theta \log(4B_\theta\Lambda_r n^{\frac12})\big)^\frac12 \\
 \leq&4n^{-\frac12}+Cn^{-\frac12}\log^2(d+2) W^{4L-4}B_\theta^{4L} N_\theta^{\frac12}\big(\log^\frac12 B_\theta+\log^\frac12 \Lambda_{r}+\log^\frac12 n\big)\\
 \leq& C n^{-\frac12} B_\theta^{4L} N_\theta^{4L-\frac{7}{2}} \big( \log^\frac12 B_\theta + \log^\frac12 N_\theta + \log^\frac12 n \big).
\end{align*}
We endow the space $\Theta \times I_{\mathbf{c}}^M\times\Omega_\gamma^M$ with
$ \|(\theta, \mathbf{c}, X)\|_*=\max\{|\theta|_{\ell_\infty}, |\mathbf{c}|_{\ell_\infty}, \max_{i,j}|X_{i,j}|\}$, for any $ (\theta, \mathbf{c}, X)\in\Theta \times I_{\mathbf{c}}^M\times\Omega_\gamma^M$. Then $\|\cdot\|_*$ is a metric and $(\Theta \times I_{\mathbf{c}}^M\times\Omega_\gamma^M, \|\cdot\|_*)$ a metric space. Then using \cite[Proposition 5]{CuckerSmale:2002} again we have
\begin{align*}
 \log \mathcal{C}(\mathcal{H}_d,\|\cdot\|_{L^{\infty}(\partial\Omega)},\epsilon)&\leq
\log \mathcal{C}(\Theta \times I_{\mathbf{c}}^M\times\Omega_\gamma^M,\|\cdot\|_* ,\Lambda_{d}^{-1}\epsilon)  \leq CN_\theta \log(4B_\theta \Lambda_d\epsilon^{-1}),
\end{align*}
with $\Lambda_{d} = CLW^{L+1}B_\theta^L$. Similarly for the third term, there holds
\begin{align*}
\log \mathcal{C}(\mathcal{H}_n,\|\cdot\|_{L^{\infty}(\partial\Omega)},\epsilon)&\leq
\log \mathcal{C}(\Theta \times I_{\mathbf{c}}^M\times\Omega_\gamma^M,\|\cdot\|_*,\Lambda_{n}^{-1}\epsilon)  \leq CN_\theta \log(4B_\theta \Lambda_n\epsilon^{-1}),
\end{align*}
with $\Lambda_n=CL^2 W^{3L-3} B_\theta^{3L-2}$. Repeating the preceding argument yields
\begin{align*}
\mathfrak{R}_n(\mathcal{H}_d)&\leq Cn^{-\frac{1}{2}}B_\theta^2 N_\theta^{\frac{5}{2}}(\log^{\frac{1}{2}} B_\theta+\log^{\frac{1}{2}} N_\theta+\log^{\frac{1}{2}} n),\\
\mathfrak{R}_n(\mathcal{H}_n)&\leq Cn^{-\frac{1}{2}}B_\theta^{2L} N_\theta^{2L-\frac{3}{2}}(\log^{\frac{1}{2}} B_\theta+\log^{\frac{1}{2}} N_\theta+\log^{\frac{1}{2}} n).
\end{align*}
Finally, the desired result follows the PAC type bound in Lemma \ref{lem:PAC}.
\end{proof}

Last, we can state the main result of this part.
\begin{theorem}\label{thm:main}
Let Assumptions \ref{assump:separation} and \ref{assump:regularity} hold. Then for any small $\epsilon,\tau,\mu>0$,  there exists an NN
$v_\theta \in \mathcal{N}_{\rho}(C, C\epsilon^{-\frac{d}{1-\mu}},C\epsilon^{-2-\frac{9d/2+4+2\mu}{1-\mu}})$, if the numbers $N_r$ and $N_b$ of sampling points satisfy $N_r,N_b\sim\mathcal{O}(\epsilon^{-\frac{4}{1-\mu}}{\log^{\frac{1}{1-\mu}}\frac{1}{\tau}}{})$, then with probability at least $1- 3 \tau$, there exists a permutation $\pi$ of the set $\{1,\ldots,M\}$ such that for $p<\frac{d}{d-1}$,
\begin{align*}	    	
\max_{1\leq j\leq M}{\rm dist}(\widehat{\mathbf{x}}_{j}^*,\mathbf{x}_{\pi(j)}^*)+
\max_{1\leq j\leq M}\!\!|\widehat{c}_j^*-c_{\pi(j)}^*|+
\|v^*-v_{\widehat{\theta}^*}\|_{L^2(\Omega)}+\|u^*-u_{\widehat{\theta}^*}\|_{L^p(\Omega)}\leq C(\epsilon+\delta)^{\tfrac{1}{M}}.
\end{align*}
\end{theorem}\begin{proof}
With the choice  $N_r,N_b\sim\mathcal{O}(\epsilon^{-\frac{4}{1-\mu}}\log^{\frac{1}{1-\mu}}\frac{1}{\tau})$, $e_r\leq C\epsilon^2$ and $e_b\leq C\epsilon^2$ holds simultaneously with probability at least $1-3\tau$, and thus the following estimate holds
$ \sup_{(v_\theta,\mathbf{c},X)\in \mathcal{W}_\epsilon\times I_\mathbf{c}^M\times \Omega_\gamma^M}|\widehat{L}^\delta(v_\theta,\mathbf{c},X)-L^\delta(v_\theta,\mathbf{c},X)|\leq C\epsilon^2$. Then combining Theorems \ref{thm:disbound}--\ref{thm:errorv}, Lemma \ref{lem:tanh-approx}, Corollary \ref{cor:erroru} with Lemma \ref{lemma:errordecom} yields the assertion.
\end{proof}

\section{Numerical experiments and discussions}\label{sec:experiment}
Now we present numerical examples to illustrate SENN. The gradient of  $v_\theta(\mathbf{x})$ with respect to the input ${\mathbf{x}}$ (i.e., spatial derivative) and that of the loss $\widehat{L}^\delta$  with respect to $(\theta,\mathbf{c},X)$  are computed via automatic differentiation using \texttt{torch.autograd}. The empirical loss $\widehat{L}^\delta(v_\theta, \mathbf{c},X)$ is minimized using Adam \cite{KingmaBa:2015}, as implemented in the SciPy library, with the default setting (tolerance =1.0e-8, no box constraint). The noisy data $f^\delta$ is generated as
$f^\delta(\mathbf{y}_j)=f(\mathbf{y}_j)(1+\delta\xi(\mathbf{y}_j))$, for each sampling point $\mathbf{y}_j\in \partial\Omega$,
where $\xi(\mathbf{y}_j)$ follows an i.i.d. uniform distribution on $[-1,1]$, and the noisy Neumann data $g^\delta$ is generated similarly. The PyTorch source code for the numerical experiments will be made available at \url{https://github.com/hhjc-web/point-source-identification}.

The first example is the standard Poisson equation with point sources.

\begin{example}\label{exam:poisson}
    Let $\Omega=(-1,1)^2$, and the harmonic function $v(x_1,x_2)=x_1^2-x_2^2$, with the information of the sources given in Table \ref{table:poisson}. Consider three cases: {\rm (i)} full Cauchy data, {\rm (ii)} Dirichlet data $f^\delta$ on the edges $\{-1,1\}\times(-1,1)$ and Neumann data $g^\delta$ on $\partial\Omega$, {\rm (iii)} Cauchy data on $\partial\Omega\backslash(-1,1)\times\{-1\}$.
\end{example}

\begin{figure}[hbt!]
	\centering\setlength{\tabcolsep}{2pt}
	\begin{tabular}{ccc}
		\includegraphics[width=.32\textwidth]  {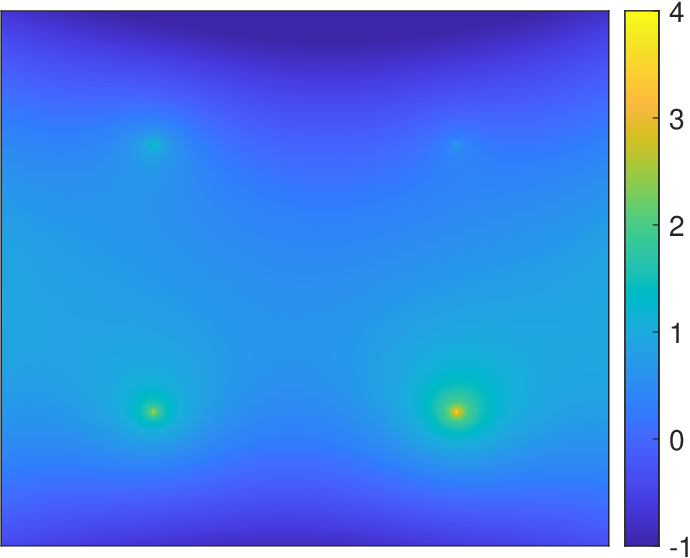} & \includegraphics[width=.32\textwidth]  {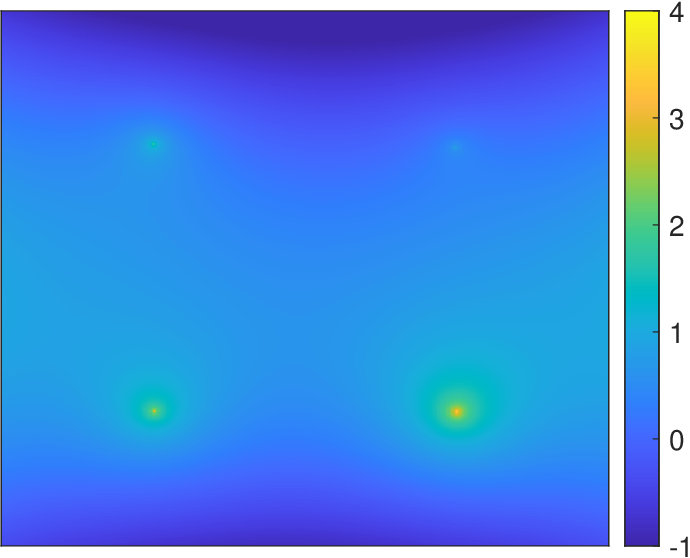} & \includegraphics[width=.33\textwidth]  {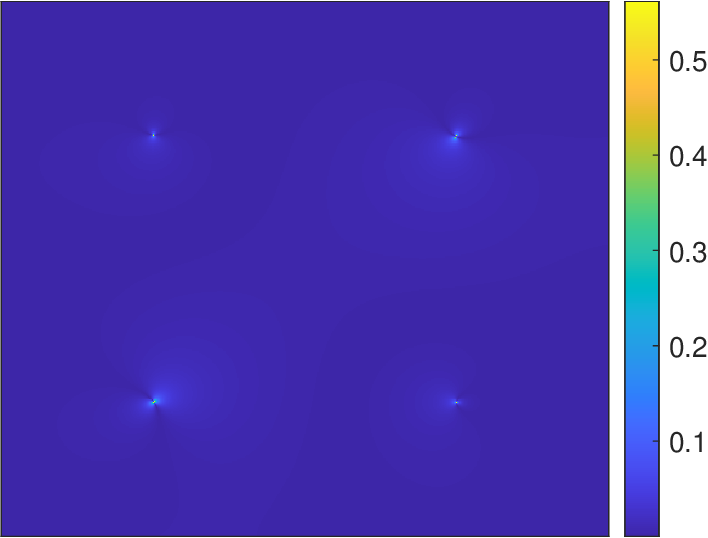}\\
		(a) exact & (b) predicted  & (c) error
	\end{tabular}
	\caption{\label{fig:exam:poisson} The neural network approximation of $u(\mathbf{x})$ for Example \ref{exam:poisson}, case (i).}
\end{figure}

First, we employ Algorithm \ref{algorithm:num} to detect the number $M$ of point sources, with a tolerance $\epsilon_{\rm tol}=0.1$ (for the noise level $\delta=2\%$). The algorithm terminates at $\widehat{M}=4$ with the singular  values $\{11.365,3.680,2.411,1.146,0.067\}$ (with 3 decimal places). We employ a 2-20-20-1 NN (2 hidden layers, each having 20 neurons) to approximate the smooth part $v(\mathbf{x})$, set the penalty weight $\sigma_d=\sigma_n=10.0$, and take $N_r=3,000$ points in the domain $\Omega$ and $N_b=2,000$ points on the boundary $\partial\Omega$ to form the empirical loss $\widehat{L}^\delta$. Further, we set the learning rate to $\text{1.0e-3}$ for NN parameters $\theta$ and $\text{6.0e-3}$ for point densities and locations. The results for the densities and locations obtained after 10,000 iterations are shown in Table \ref{table:poisson}, which also includes the results by the algebraic method  \cite{el2000inverse}. Both the proposed method and the algebraic method can give satisfactory results. The recovered $u(\mathbf{x})$ in Fig. \ref{fig:exam:poisson} shows that the pointwise error is small and the predicted locations of the point sources are quite accurate.

\begin{table}[hbt!]
\centering
\begin{threeparttable}
\caption{The exact and estimated locations and densities of point sources for Example \ref{exam:poisson} with $\delta=2\%$. Row (i) is for full Cauchy data, and rows (ii) and (iii) are for partial Cauchy data.}\label{table:poisson}
\centering
\begin{tabular}{c|cc|cc}
\toprule
point & density & location & density & location \\
\midrule
exact & 1.000 & $(0.500,0.500)$ & 2.000 & $(-0.500,0.500)$ \\
\hline
(i), SENN & 1.025 & $(0.496,0.491)$ & 1.983 & $(-0.500,0.503)$ \\
\hline
(i), \cite{el2000inverse} & 0.977 & $(0.506,0.507)$ & 2.047 & $(-0.502,0.490)$ \\
\hline
(ii) & 1.096 & $(0.458,0.478)$ & 1.962 & $(-0.497,0.504)$ \\
\hline
(iii) & 1.026 & $(0.495,0.494)$ & 1.975 & $(-0.503,0.503)$ \\
\midrule
point & density & location & density & location \\
\midrule
exact & 3.000 & $(-0.500,-0.500)$ & 4.000 & $(0.500,-0.500)$ \\
\hline
(i), SENN & 3.022 & $(-0.497,-0.498)$ & 3.973 & $(0.501,-0.500)$ \\
\hline
(i), \cite{el2000inverse} & 2.989 & $(-0.506,-0.496)$ & 3.990 & $(0.502,-0.498)$ \\
\hline
(ii) & 3.038 & $(-0.499,-0.498)$ & 3.899 & $(0.509,-0.501)$ \\
\hline
(iii) & 3.126 & $(-0.487,-0.503)$ & 3.844 & $(0.511,-0.496)$ \\
\bottomrule
\end{tabular}
\end{threeparttable}
\end{table}

The smallest singular value $\sigma_s(A)$ of the Hankel matrix $A$ in Algorithm \ref{algorithm:num} increases with the noise level $\delta$. If the noise level $\delta$ is very high, or the tolerance $\epsilon_{\rm tol}$ is not properly set, the number $M$ of point sources may be inaccurately estimated, which can affect the accuracy of the subsequent recovery. Nonetheless, Theorem \ref{thm:main} shows that the  Hausdorff distance $d_H(\widehat{X}^*, X^*)$ can still be properly bounded (under suitable conditions). Numerically, we test the noise level $\delta=10\%$ and choose different $\epsilon_{\rm tol}$ so that the number $\widehat{M}$ of point sources is incorrectly estimated to be 3 and 5 respectively, and then proceed to training. The estimated locations and densities of the point sources are shown in Table \ref{table:poisson-wrong}. The results indicate that for $\widehat{M} >M $, the predicted point sources are still very satisfactory: it roughly amounts to splitting one or more sources, and their combined effect is close to the exact one. However, for $\widehat{M}<M$, the predictions are not so satisfactory. Thus, when the noise level $\delta$ is high, one should set the tolerance $\epsilon_{\rm tol}$ slightly smaller, so that the number $\widehat{M}$ of predicted points is no less than the exact one $M$. We repeat the experiment with the algebraic method \cite{el2000inverse}; see Table \ref{table:poisson-wrong2} for the results. It is observed that when the number $M$ of points is incorrectly estimated, the error is very large and the prediction is not reliable. This is attributed to the ill-conditioning of the Hankel matrix $A$ used in the algebraic method.

\begin{table}[hbt!]
\centering
\begin{threeparttable}
\caption{\label{table:poisson-wrong} The exact and estimated locations and densities of point sources for Example \ref{exam:poisson}, case (i) with $\delta=10\%$, using SENN.}
\centering
\setlength{\tabcolsep}{5pt}
\begin{tabular}{cc|cc|cc}
\toprule
\multicolumn{2}{c}{exact} & \multicolumn{2}{c}{$\widehat M= 5$}& \multicolumn{2}{c}{$\widehat M= 3$}\\

\midrule
density & location & density & location & density & location \\
\midrule
$\begin{matrix}
    1.000\\2.000\\3.000\\4.000
\end{matrix}$ &
$\begin{matrix}
    (0.500,0.500)\\(-0.500,0.500)\\(-0.500,-0.500)\\(0.500,-0.500)
\end{matrix}$ &
$\begin{matrix}
    1.054\\2.051\\3.104\\1.826\\2.051
\end{matrix}$ &
$\begin{matrix}
    (0.499,0.476)\\(0.460,-0.542)\\(-0.485,-0.484)\\(0.561,-0.471)\\(0.460,-0.542)
\end{matrix}$ &
$\begin{matrix}
    3.906\\2.142\\3.831
\end{matrix}$ &
$\begin{matrix}
    (-0.221,0.252)\\(-0.523,-0.542)\\(0.501,-0.493)
\end{matrix}$ \\
\bottomrule
\end{tabular}
\end{threeparttable}
\end{table}

\begin{table}[hbt!]
\centering
\begin{threeparttable}
\caption{\label{table:poisson-wrong2} The exact and estimated locations and densities of point sources for Example \ref{exam:poisson}, case (i) with $\delta=10\%$, using the direct method \cite{el2000inverse}.}
\centering
\begin{tabular}{cc|cc}
\toprule
\multicolumn{2}{c}{exact} & \multicolumn{2}{c}{$\widehat M= 5$}\\

\midrule
density & location & density & location \\
\midrule
$\begin{matrix}
    1.000\\2.000\\3.000\\4.000
\end{matrix}$ &
$\begin{matrix}
    (0.500,0.500)\\(-0.500,0.500)\\(-0.500,-0.500)\\(0.500,-0.500)
\end{matrix}$ &
$\begin{matrix}
    0.628\\0.202\\2.467\\3.626\\3.086
\end{matrix}$ &
$\begin{matrix}
    (0.524,0.627)\\(-0.776,0.655)\\(-0.571,-0.482)\\(0.522,-0.501)\\(-0.385,0.309)
\end{matrix}$
\\
\bottomrule
\end{tabular}
\end{threeparttable}
\end{table}

\begin{figure}[hbt!]
	\centering\setlength{\tabcolsep}{0pt}
	\begin{tabular}{cc}
		\includegraphics[width=.45\textwidth]  {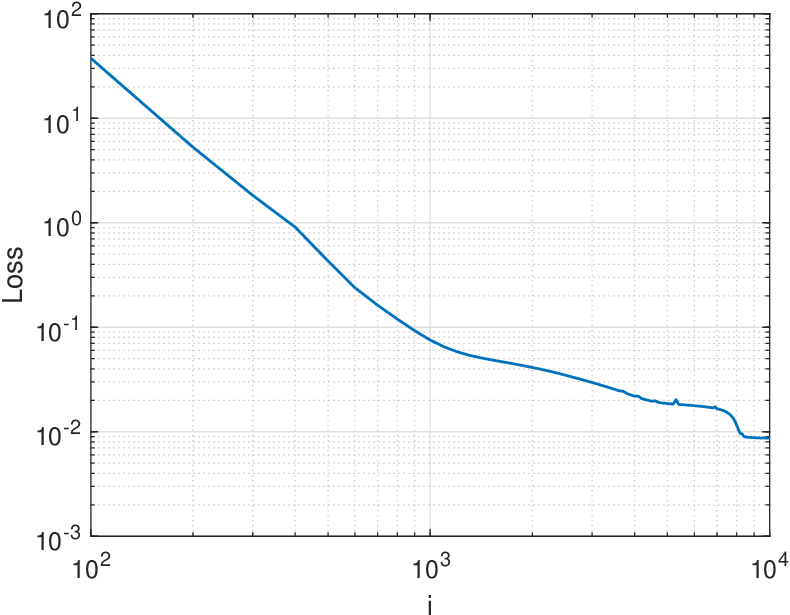} &  \includegraphics[width=.44\textwidth]  {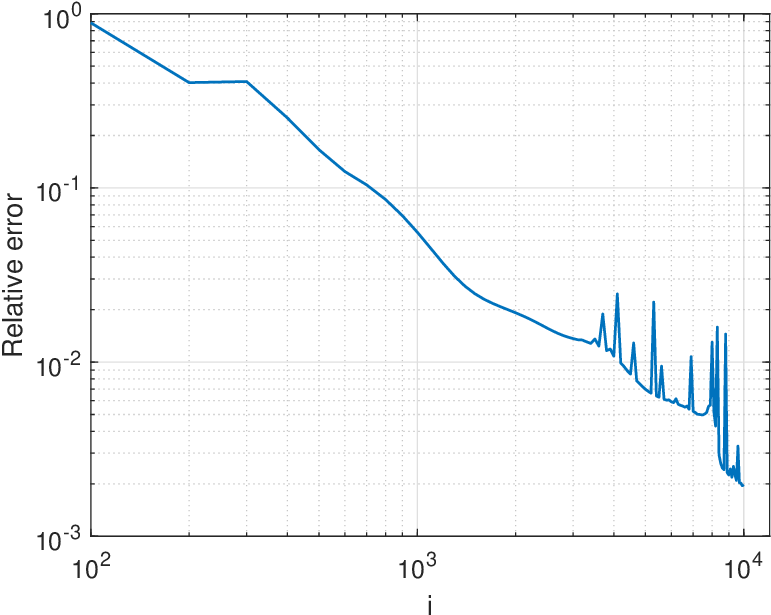}\\
		(a) loss $\widehat{L}^\delta$ vs $i$ & (b) relative error $e$ vs $i$
	\end{tabular}
	\caption{\label{fig:exam:poisson:decay} The dynamics of the training process for the proposed method for Example \ref{exam:poisson}, case (i): {\rm(a)} the decay of the loss $\widehat{L}^\delta$ versus the iteration index $i$, {\rm(b)} the error $e$ versus the iteration index $i$.}
\end{figure}

\begin{figure}[hbt!]
	\centering\setlength{\tabcolsep}{2pt}
	\begin{tabular}{cccc}
	\includegraphics[height=4.5cm]  {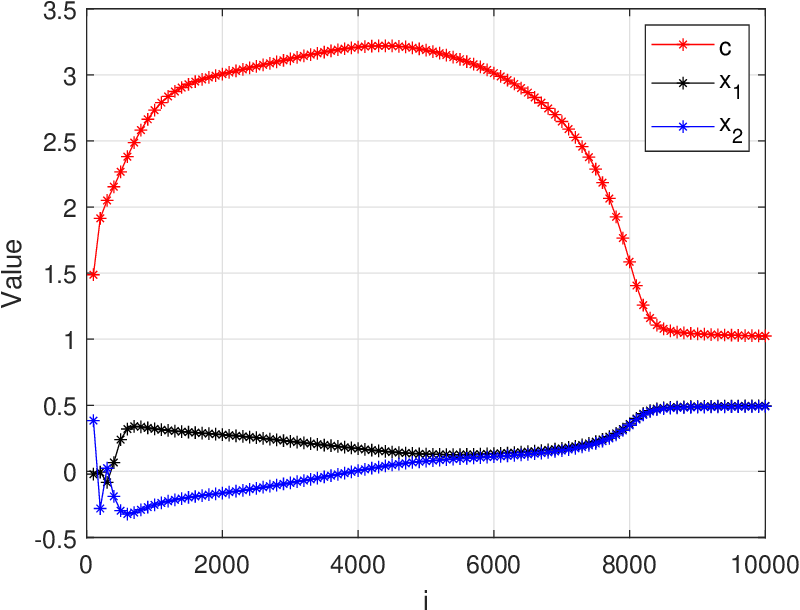} &
      \includegraphics[height=4.5cm]  {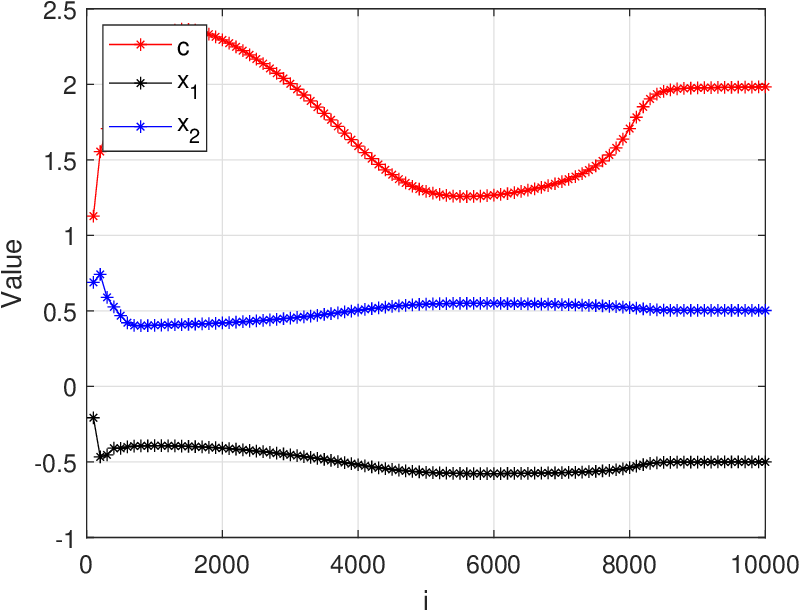}  \\
       (a) point 1 & (b) point 2 \\
        \includegraphics[height=4.5cm]  {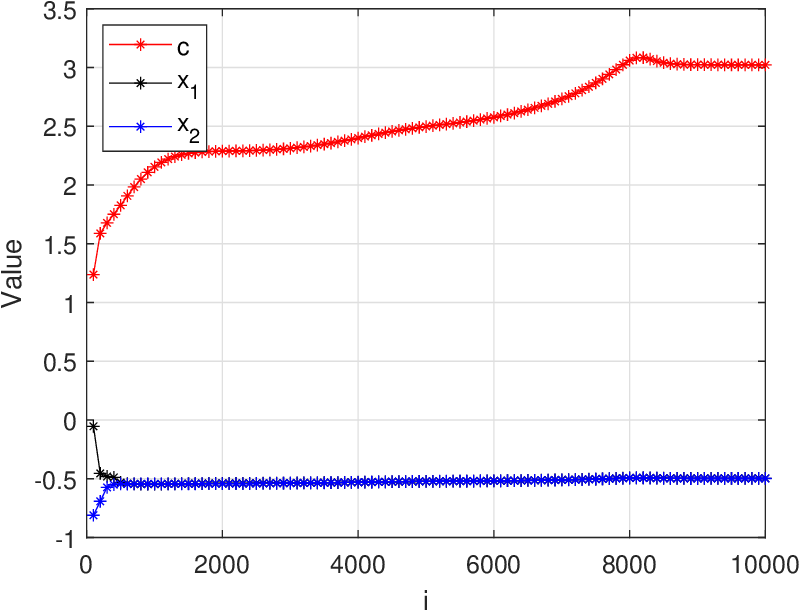} &
       \includegraphics[height=4.5cm]  {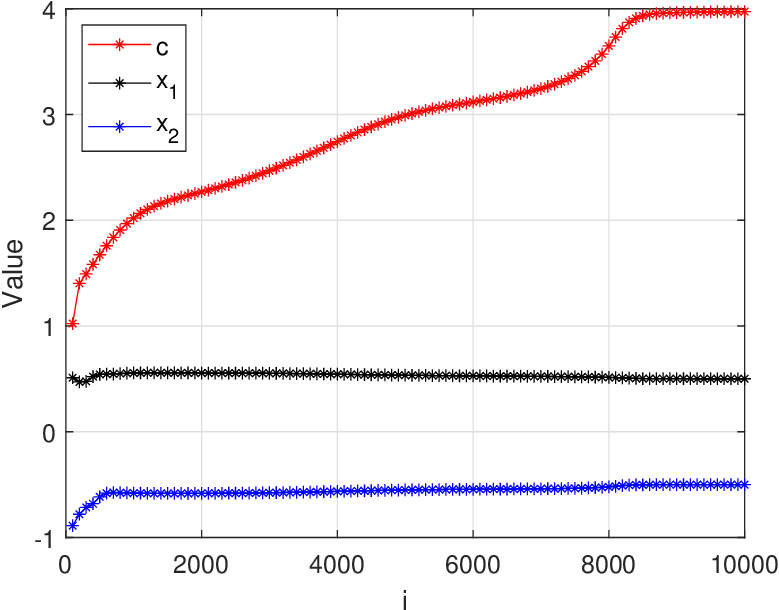} \\
     (c) point 3 & (d) point 4
	\end{tabular}
	\caption{\label{fig:exam:poisson:dynamic} The training dynamics of the density $c$ and location $\mathbf{x}=(x_1,x_2)$ of the four point sources for Example \ref{exam:poisson}, case (i).}
\end{figure}

To shed further insights, we show in Fig. \ref{fig:exam:poisson:decay} the training dynamics of the empirical loss $\widehat{L}^\delta$ and the relative error $e=\|v_\theta-v^*\|_{L^2(\Omega)}/\|v^*\|_{L^2(\Omega)}$ (of the regular part $v_\theta$), where $i$ denotes the total iteration index. Both the loss $\widehat{L}^\delta$ and the relative error $e$ decay steadily as the iteration proceeds, indicating stable convergence of the Adam optimizer. However, as the number of iterations increases further, the loss trajectory exhibits oscillations while the error $e$ stays at a low level. This is typical for gradient type optimizers such as  Adam \cite{BottouCurtis:2018}: Large learning rates or gradients can cause oscillations. The final error $e$ saturates at around $ 10^{-3} $, which agrees with the observation that the accuracy of neural PDE solvers tends to stagnate at a level of $10^{-2}\sim 10^{-3}$ \cite{WangPerdikaris:2022jcp}.
Fig. \ref{fig:exam:poisson:dynamic} shows the dynamics of the recovered density $c_j$ and location $\mathbf{x}_j$ of each point source during the training, where $\mathbf{x}_j$ is initialized randomly in the domain $\Omega$, and the density $c_j$ is initialized randomly in $(0,1)$. Numerically, we observe that the densities $c_j$ takes longer time to reach convergence than the locations $\mathbf{x}_j$, but after 10,000 iterations both densities and locations stabilize at the minimizer.

In practice, we may only know partial Cauchy data, for which we illustrate two cases (ii) and (iii). Since the algebraic method does not work any more for partial data, we assume that the number $M$ of point sources is given. To form the empirical loss $\widehat{L}^\delta$, for case (ii), we take $1,000$ observation points for Dirichlet data, $2,000$ for Neumann data, and $N_r=5,000$ in the domain $\Omega$, and for case (iii), we take the number of observation points $N_b=9,000$ on the boundary and $N_r=10,000$ in the domain $\Omega$. The recovery results in Table \ref{table:poisson} (partial data 1) show that the method works reasonably well for partial data, but the results are slightly worse than case (i).

The next example is about a high-dimensional problem.
\begin{example}\label{exam:poisson10}
Consider Example \eqref{problem} with $\Omega=(-1,1)^{10}$. Let $v(\mathbf{x})=-9x_1^2+\sum_{i=2}^{10}x_i^2$, and the exact point sources are shown in Table \ref{table:poisson10}.
\end{example}

In the experiment, assuming a known $M$, we directly train the NN. To study the impact of data noise on the  recovery accuracy, we test several noise levels. We employ an NN architecture 10-50-50-50-50-50-50-1 to approximate the regular part $v(\mathbf{x})$, the penalty factor $\sigma_d=\sigma_n=20.0$, and take $N_r=10,000$ points in the domain $\Omega$ and $N_b=20,000$ on the boundary $\partial\Omega$ to form the empirical loss $\widehat{L}^\delta$. The learning rate for NN parameters is 1e-3, and that for locations and densities varies with the noise level $\delta$, cf. Table \ref{table:poisson10}. When the noise level $\delta$ is large, the loss $\widehat L^\delta$ has a large minimum value, but the relative error $e$ of the approximation $v_{\hat\theta^*}$ is still relatively small. Numerically, we observe that the inaccuracy of the approximation $v_{\widehat \theta^*}$concentrates on the boundary $\partial\Omega$, and the internal accuracy is relatively high. The results in Table \ref{table:poisson10} agree with the intuition: the larger $\delta$ is, the less accurate the recovered regular part $v_{\widehat\theta^*}$ is. The recovery accuracy for exact data can reach nearly $10^{-3}$, showing the excellent capability of the proposed SENN for recovering point sources.

\begin{table}[hbt!]
\centering\setlength{\tabcolsep}{5pt}
\begin{threeparttable}
\caption{\label{table:poisson10} The exact and estimated locations and densities of point sources for Example \ref{exam:poisson10} using SENN. The learning rate for the locations and intensities when $\delta=0\%$, 2\%, 5\% and 10\% is 6e-3, 2e-3, 8e-3 and 2e-3, respectively.}
\centering
\begin{tabular}{c|c|c|c|c|c}
\toprule
\multirow{2}{*}{point} & \multirow{2}{*}{exact} & \multicolumn{4}{c}{predicted} \\
\cline{3-6}
&& $\delta=0\%$ & $\delta=2\%$ & $\delta=5\%$ & $\delta=10\%$\\
\midrule
$c_1$ & 100.0 & 98.9 & 93.4 & 87.5 & 80.5 \\
\hline
$\mathbf{x}_1$ & $\begin{pmatrix}
    -0.500\\-0.500\\0.000\\0.000\\0.000\\0.000\\0.000\\0.000\\0.000\\0.000
\end{pmatrix}$ & $\begin{pmatrix}
    -0.505\\-0.500\\0.002\\0.001\\0.001\\-0.002\\-0.002\\0.002\\0.003\\-0.001
\end{pmatrix}$ & $\begin{pmatrix}
    -0.503\\-0.500\\0.024\\-0.002\\0.016\\-0.006\\-0.015\\-0.004\\-0.002\\0.003
\end{pmatrix}$ & $\begin{pmatrix}
    -0.505\\-0.491\\0.054\\-0.066\\0.041\\-0.019\\-0.034\\-0.023\\0.000\\-0.001
\end{pmatrix}$ & $\begin{pmatrix}
    -0.571\\-0.381\\0.153\\-0.081\\0.123\\-0.071\\-0.058\\-0.083\\-0.048\\-0.065
\end{pmatrix}$\\
\midrule
$c_2$ & 100.0 & 100.6 & 96.4 & 93.5 & 77.7 \\
\hline
$\mathbf{x}_2$ & $\begin{pmatrix}
    0.500\\0.500\\0.000\\0.000\\0.000\\0.000\\0.000\\0.000\\0.000\\0.000
\end{pmatrix}$ & $\begin{pmatrix}
    0.497\\0.503\\-0.001\\-0.003\\0.003\\-0.002\\0.001\\-0.002\\-0.001\\0.000
\end{pmatrix}$ & $\begin{pmatrix}
    0.482\\0.502\\0.015\\-0.014\\0.009\\0.014\\0.024\\-0.025\\-0.007\\-0.013
\end{pmatrix}$ & $\begin{pmatrix}
    0.470\\0.491\\0.053\\-0.033\\0.018\\0.036\\0.059\\-0.062\\-0.021\\-0.041
\end{pmatrix}$ & $\begin{pmatrix}
    0.499\\0.463\\0.118\\-0.111\\0.030\\0.043\\0.099\\-0.129\\-0.038\\-0.090
\end{pmatrix}$\\
\bottomrule
\end{tabular}

\end{threeparttable}

\end{table}

Last, we extend the method to the Helmholtz equation.
\begin{example}\label{exam:hel}
Consider the inverse source problem with the Helmholtz equation
$\Delta u+k^2u=\sum_{j=1}^3c_j\delta_{\mathbf{x}_j}$, for $\mathbf{x}_1,\mathbf{x}_2,\mathbf{x}_3\in\Omega=:( -1,1)^3\subset\mathbb{R}^3$,
with noisy Cauchy data. The fundamental solution $\Phi(\mathbf{x})$ is given by
$\Phi(\mathbf{x})=\frac{{\rm e}^{{\rm i}k|\mathbf{x}|}}{4\pi|\mathbf{x}|}$. We set the frequency $k=1$ and regular part $v(\mathbf{x})=\sin(\tfrac{x_1+x_2+x_3}{\sqrt{3}})$. The point sources, i.e., $(c_j,\mathbf{x}_j)$, are shown in Table \ref{table:exam:hel}. Consider two cases: {\rm(i)} full Cauchy data, {\rm(ii)} Cauchy data on $\partial\Omega\setminus (-1,1)^2\times\{-1\}$.
\end{example}

\begin{figure}[hbt!]
	\centering\setlength{\tabcolsep}{2pt}
	\begin{tabular}{ccc}
		\includegraphics[width=.3\textwidth]  {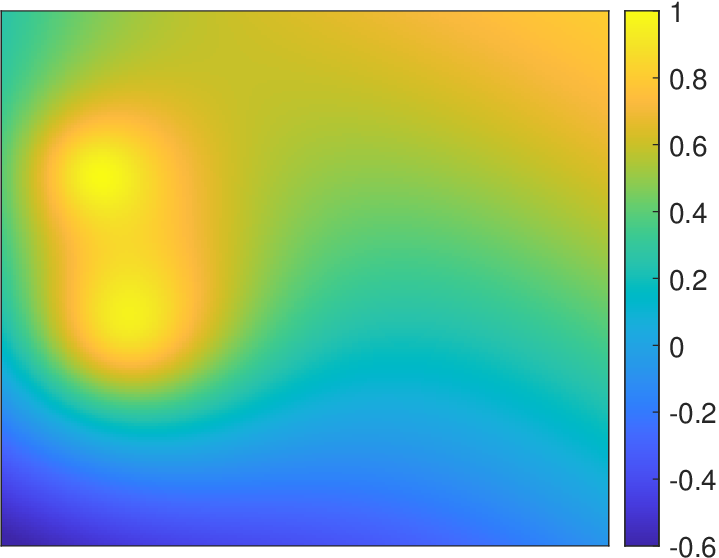} & \includegraphics[width=.3\textwidth]  {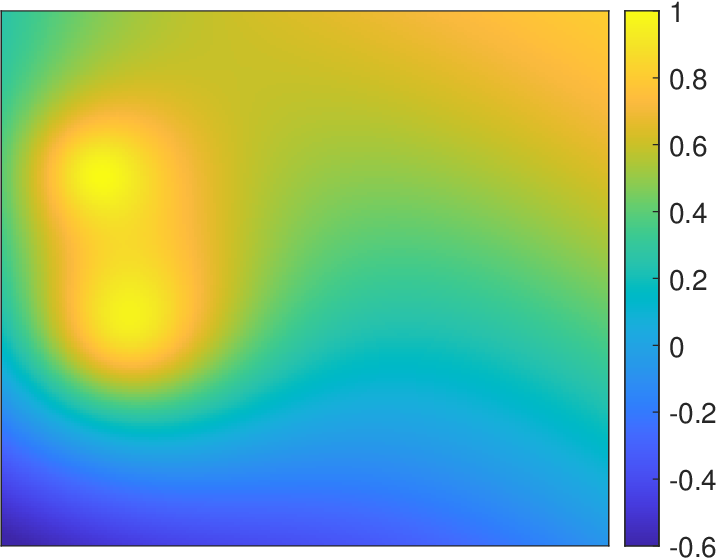} & \includegraphics[width=.31\textwidth]  {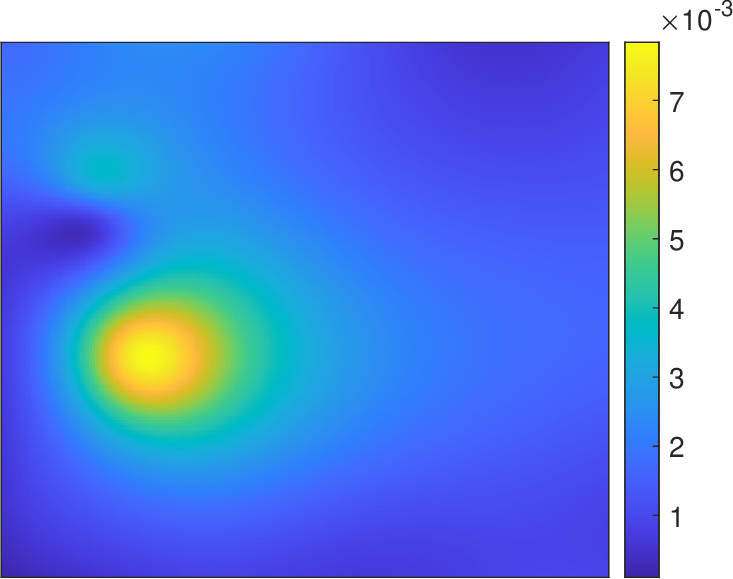}\\
        \includegraphics[width=.3\textwidth]  {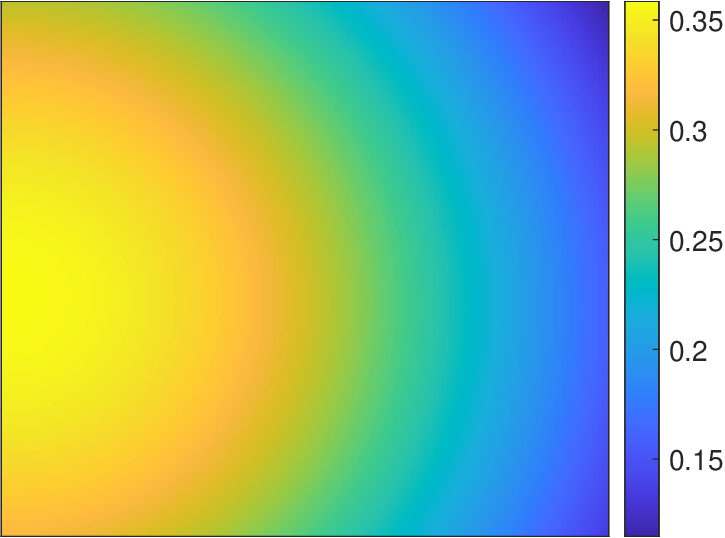} & \includegraphics[width=.3\textwidth]  {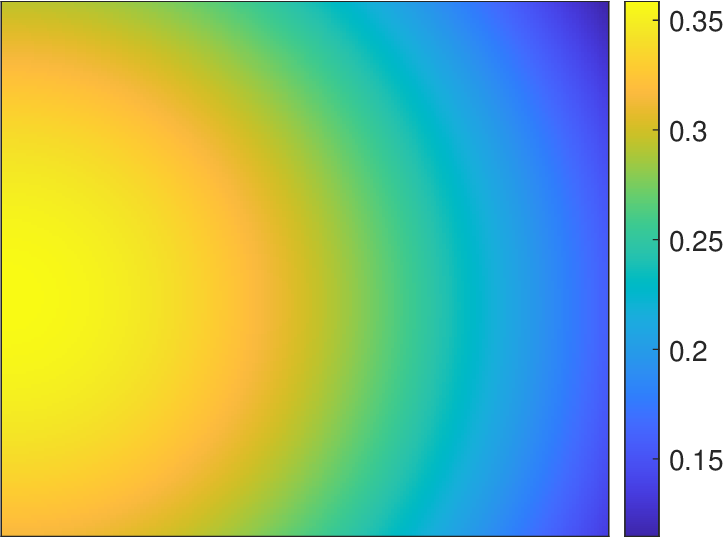} & \includegraphics[width=.3\textwidth]  {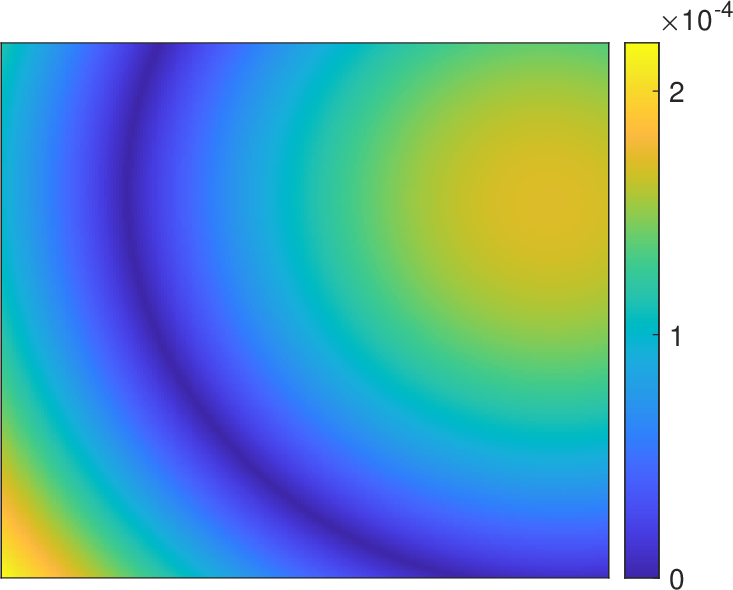}\\
		(a) exact & (b) predicted  & (c) error
	\end{tabular}
	\caption{\label{fig:exam:hel} The NN approximation of the real (top) and imaginary (bottom) parts of $u(\mathbf{x})$ for Example \ref{exam:hel}, case (i), at the slice $x_3=0$.}
\end{figure}

The difference between Helmholtz and Poisson cases lies in the fact that the former is complex-valued. We again assume a known number $M$ of point sources, and employ $N_r=10,000$ points in $\Omega$ and $N_b=20,000$ points on $\partial\Omega$ to form the empirical loss $\widehat{L}^\delta$. We use two 3-20-20-20-1 NNs to approximate the real and imaginary parts of the regular part $v(\mathbf{x})$, take the penalty weights $\sigma_d=\sigma_n=20.0$, and set the learning rate to $\text{1.0e-3}$ for $\theta$ and $\text{2.0e-3}$ for $(\mathbf{c},X)$. The training results after 10,000 iterations are shown in Table \ref{table:exam:hel} and Fig. \ref{fig:exam:hel} for the noise level $\delta=5\%$.  Fig. \ref{fig:exam:hel} shows the real and imaginary parts of $u(\mathbf{x})$ sliced at $x_3=0$. The pointwise error is small. For partial Cauchy data, the recovered intensities and locations are still satisfactory, cf. Table \ref{table:exam:hel}. This again shows the capability of SENN for identifying point sources.

\begin{table}[hbt!]
\centering
\begin{threeparttable}
\caption{\label{table:exam:hel} The exact and estimated locations and densities of point sources for Example \ref{exam:hel} (with noise level $\delta=5\%$). Rows (i) and (ii) denote full and partial Cauchy data.}
\centering
\begin{tabular}[5pt]{c|cc|cc|cc}
\toprule
point & $c_1$ &$\mathbf{x}_1$  & $c_2$ &$\mathbf{x}_2$ & $c_3$ &$\mathbf{x}_3$ \\
\midrule
exact & 5.000 & $\begin{pmatrix}
    -0.600\\-0.200\\0.300
\end{pmatrix}$ & 2.000 & $\begin{pmatrix}
    -0.700\\0.400\\-0.200
\end{pmatrix}$ & -3.000 &$\begin{pmatrix}
    0.000\\0.100\\0.900
\end{pmatrix}$\\
\hline
(i) & 5.016 & $\begin{pmatrix}
    -0.599\\-0.200\\0.299
\end{pmatrix}$ & 1.999 & $\begin{pmatrix}
    -0.700\\0.401\\-0.200
\end{pmatrix}$ & -3.022 & $\begin{pmatrix}
    0.000\\0.100\\0.899
\end{pmatrix}$\\
\hline
(ii) & 5.009 & $\begin{pmatrix}
    -0.600\\-0.200\\0.300
\end{pmatrix}$ & 2.006 & $\begin{pmatrix}
    -0.700\\0.399\\-0.199
\end{pmatrix}$ & -3.026&  $\begin{pmatrix}
    0.000\\0.101\\0.899
\end{pmatrix}$\\
\bottomrule
\end{tabular}
\end{threeparttable}
\end{table}

\begin{figure}[hbt!]
	\centering\setlength{\tabcolsep}{0pt}
	\begin{tabular}{ccc}
		\includegraphics[height=3.5cm]  {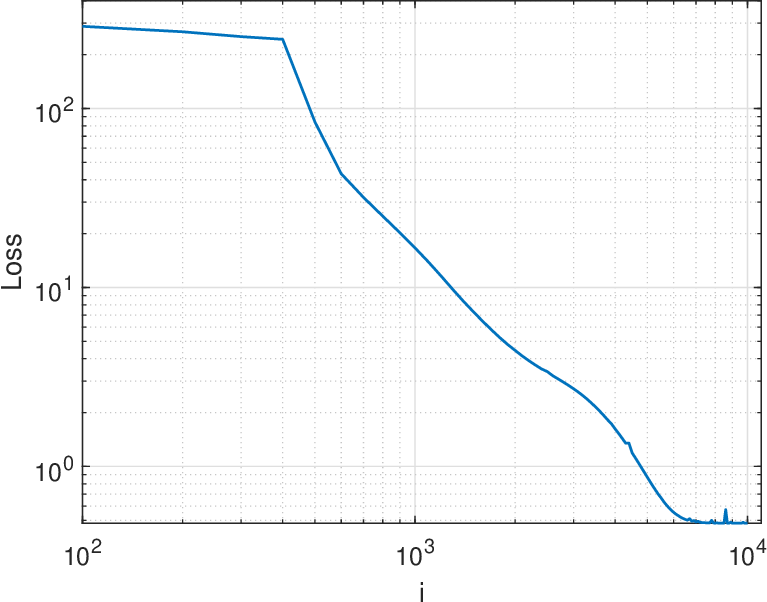} &
        \includegraphics[height=3.5cm]  {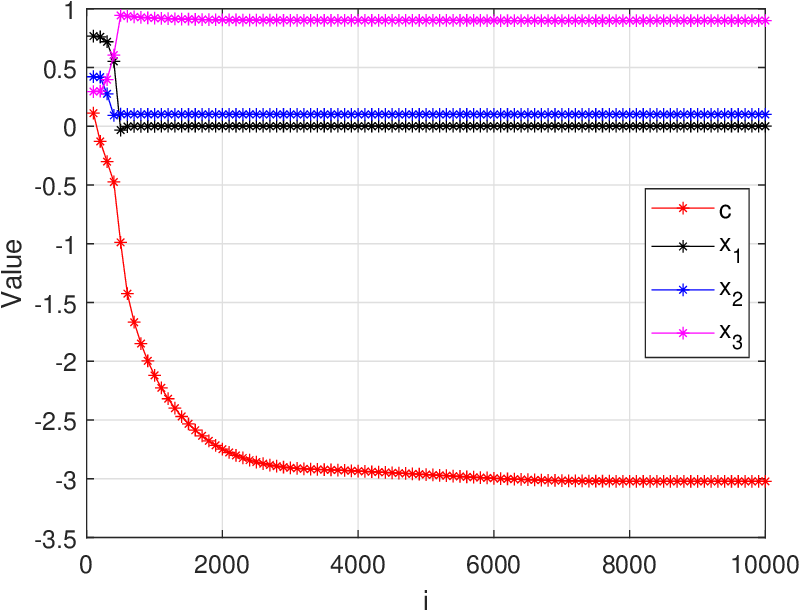} &
        \includegraphics[height=3.5cm]  {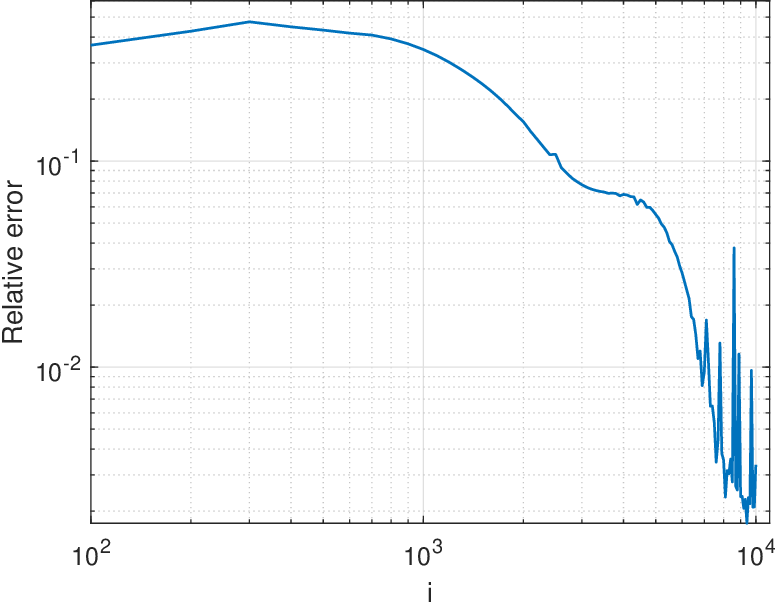}\\
		(a) $\widehat{L}^\delta$ vs $i$ & (b) point 3 & (c) $e$ vs $i$
	\end{tabular}
	\caption{\label{fig:exam:hel:decay} The dynamics of the training process for the proposed method for Example \ref{exam:hel}, case (i): {\rm(a)} the decay of the loss $\widehat{L}^\delta $ versus the iteration index $i$; {\rm(b)} the training dynamics of one singular point (i.e. the intensity $c$ and the locations) ; {\rm(c)} the relative error $e$ versus the iteration index $i$.}
\end{figure}

Fig. \ref{fig:exam:hel:decay} shows the dynamics of the loss $\widehat{L}^\delta$ and the relative error $e$ during the training process. We observe from Figs. \ref{fig:exam:hel:decay}(a) that the decay of the loss $\widehat{L}^\delta$ accelerates after around 400 iterations. This coincides with Fig. \ref{fig:exam:hel:decay}(b), where the location of one of the points converges, as indicated by the intersection of the iteration trajectory.
Like in Example \ref{exam:poisson}, after 5,000 iterations, the iteration still exhibits oscillations in both loss curve and relative error, akin to the catapult phenomenon in deep learning.

Last, we give one example with a variable coefficient.
\begin{example}\label{exam:variable}
Consider the following second-order elliptic equation with a variable coefficient: $- \nabla \cdot(\kappa\nabla u) = \sum_{j=1}^2c_j\delta_{\mathbf{x}_j}(\mathbf{x})$ in $ \Omega=:(-1,1)^2$ with $\kappa=x_1^2+x_2^2+1$. The information of the source is shown in Table \ref{table:variable}.
\end{example}

\begin{table}[hbt!]
\centering
\begin{threeparttable}
\caption{\label{table:variable} The exact and estimated locations and densities of point sources for Example \ref{exam:variable} with different noise level $\delta$.}
\centering
\begin{tabular}{c|cc|cc}
\toprule
point & density & location & density & location \\
\midrule
exact & 10.000 & $(0.333,0.333)$ & 10.000 & $(-0.333,-0.333)$ \\
\hline
$\delta=0\%$ & 9.726 & $(0.347,0.344)$ & 10.172 & $(-0.320,-0.332)$ \\
\hline
$\delta=2\%$ & 9.717 & $(0.347,0.343)$ & 10.240 & $(-0.316,-0.327)$ \\
\hline
$\delta=5\%$ & 9.703 & $(0.3448,0.343)$ & 10.287 & $(-0.319,-0.323)$ \\
\bottomrule
\end{tabular}
\end{threeparttable}
\end{table}

We equip the equation with a Neumann boundary condition $\kappa\frac{\partial u}{\partial n} = -\frac{1}{|\partial\Omega|}\sum_{j=1}^2c_j$ on $ \partial\Omega$, and fix $u(1,1)=0$ in order to ensure the uniqueness of the direct problem.
Since the exact solution $u$ is unavailable, we solve the direct problem using the Galerkin finite element method (FEM) first and then use the Dirichlet data to detect the sources.
In order to get accurate data, we use the singularity splitting strategy \cite{gjerde2019splitting, HuJinZhou:2022} by splitting the solution $u$ into
$u=\sum_{j=1}^2\frac{c_j}{\kappa}\Phi_{\mathbf{x}_j}+v$
and approximating $v$ using the FEM. Due to the presence of the variable coefficient $\kappa(\mathbf{x})$, the direct method in \cite{el2000inverse} no longer works. The numerical results in Table \ref{table:variable} show that the method works equally well for the variable coefficient $\kappa(\mathbf{x})$, which shows its versatility for more general problems.

\section{Conclusions}\label{sec:concl}
In this work, we have developed a novel neural network based algorithm to identify point sources in the Poisson equation from Cauchy data. It employs the fundamental solution of the elliptic operator to extract the singular part of the solution and neural networks to approximate the regular part, and to estimate the solution (state) and unknown parameters of the point sources by optimizing an empirical loss. We have also derived an error estimate for the method. This is achieved by first establishing H\"{o}lder conditional stability for the recovered sources and densities, and the solution of the differential equation in terms of the loss function, and then deriving a priori bound on the empirical loss using techniques in statistical learning theory. We presented several numerical experiments to show its performance, which show satisfactory recovery for point sources and excellent stability with respect to noise.

There are several directions for further research. First, it is of interest to extend the approach to more general sources, e.g., that concentrate on lines, curves or surfaces, which is prevalent in  biological modeling, and to develop relevant theoretical analysis, including partial boundary data. Second, it is important to analyze other types of sources, e.g., dipoles or mixture of monopoles and dipoles, or more general distributed sources. Last, it is important to study other type of mathematical models, e.g., parabolic or hyperbolic systems, which correspond to the nonstationary counterpart of the elliptic case and potentially involve time-dependent intensities and moving locations.

\appendix

\section{Technical estimates}\label{app:est}

To bound the terms in the splitting \eqref{eqn:err-sta}, we employ Rademacher complexity of the function classes $\mathcal{H}_r$, $\mathcal{H}_d$ and $\mathcal{H}_n$. Rademacher complexity \cite{AnthonyBartlett:1999,BartlettMendelson:2002} measures the complexity of a collection of functions
by the correlation between function values with Rademacher random variables, i.e., with probability $P(\omega=1)=P(\omega=-1)=\frac12$.
\begin{definition}\label{def: Rademacher}
Let $\mathcal{F}$ be a real-valued function class defined on the domain $D$ and $\xi=\{\xi_j\}_{j=1}^n$ be i.i.d. samples from a distribution $\mathcal{U}(D)$.
Then the Rademacher complexity $\mathfrak{R}_n(\mathcal{F})$ of $\mathcal{F}$ is defined by
\begin{equation*}
\mathfrak{R}_n(\mathcal{F})=\mathbb{E}_{\xi,\omega}\bigg{[}\sup_{f\in\mathcal{F}}\ n^{-1}\bigg{\lvert}\ \sum_{j=1}^{n}\omega_j f(\xi_j)\ \bigg{\rvert} \bigg{]},
\end{equation*}
where $\omega=\{\omega_j\}_{j=1}^n$ are i.i.d Rademacher random variables.
\end{definition}

We have the following PAC-type generalization bound \cite[Theorem 3.1]{Mohri:2018}.
\begin{lemma}\label{lem:PAC}
Let $X_1,\ldots,X_n$ be a set of i.i.d. random variables, and let $\mathcal{F}$
be a function class defined on $D$ such that $\sup_{f\in\mathcal{F}}\|f\|_{L^\infty(D)}\leq M_\mathcal{F}<\infty$. Then for any $\tau\in(0,1)$, with probability at least
$1-\tau$:
\begin{equation*}
  \sup_{f\in \mathcal{F}}\bigg|n^{-1}\sum_{j=1}^n f(X_j)-\mathbb{E}[f(X)]\bigg| \leq 2\mathfrak{R}_n(\mathcal{F}) + 2M_\mathcal{F}\sqrt{\frac{\log\frac{1}{\tau}}{2n}}.
\end{equation*}
\end{lemma}

Next, we state Dudley's lemma (\cite[Theorem 9]{lu2021priori}, \cite[Theorem 1.19]{wolf2018mathematical}), which bounds Rademacher complexities using covering number.
\begin{definition}\label{def:cover}
Let $(\mathcal{M},m)$ be a metric space of real valued functions, and $\mathcal{G}\subset \mathcal{M}$. A set $\{x_i\}_{i=1}^n \subset \mathcal{G}$ is called an $\epsilon$-cover of $\mathcal{G}$ if for any $x\in \mathcal{G}$, there exists an element $x_i\in \{x_i\}_{i=1}^n$ such that $m(x, x_i) \leq \epsilon$. The
$\epsilon$-covering number $\mathcal{C}(\mathcal{G}, m, \epsilon)$ is the minimum
cardinality among all $\epsilon$-covers of $\mathcal{G}$ with respect to $m$.
\end{definition}

\begin{lemma}\label{lem:Dudley}
Let $\mathcal{F}$ be a real-valued function class defined on the domain $D$, $M_\mathcal{F}:=\sup_{f\in\mathcal{F}} \|f\|_{L^{\infty}(D)}$, and $\mathcal{C}(\mathcal{F},\|\cdot\|_{L^{\infty}(D)},\epsilon)$ the covering number of $\mathcal{F}$. Then  $\mathfrak{R}_n(\mathcal{F})$ is bounded by
\begin{equation*}
\mathfrak{R}_n(\mathcal{F})\leq\inf_{0<s< M_\mathcal{F}}\bigg(4s\ +\ 12n^{-\frac12}\int^{M_\mathcal{F}}_{s}\big(\log \mathcal{C}(\mathcal{F},\|\cdot\|_{L^{\infty}(D)},\epsilon)\big)^{\frac12}\ {\rm d}\epsilon\bigg).
\end{equation*}
\end{lemma}

To apply Lemma \ref{lem:PAC}, we bound the Rademacher complexities $\mathcal{R}_n(\mathcal{H}_r)$, $\mathcal{R}_n(\mathcal{H}_d)$ and $\mathcal{R}_n(\mathcal{H}_n)$.
This follows from Dudley's formula in Lemma \ref{lem:Dudley}.
The next lemma gives useful boundedness and Lipschitz continuity of the NN function class in terms of $\theta$; see \cite[Lemma 3.4 and Remark 3.3]{JinLiLu:2022} and \cite[Lemma 5.3]{Jin:DNN-Control}.
These estimates also hold when $L^\infty(\Omega)$ is replaced with $L^\infty(\partial\Omega)$. The condition $\Omega\subset (-1,1)^d$ is used in deriving the bounds in Lemma \ref{lem:NN-Lip}.

\begin{lemma}\label{lem:NN-Lip}
Let $L$, $W$ and $B_\theta$ be the depth, width and maximum weight bound of an NN function class $\mathcal{W}$, with $N_\theta$ nonzero weights. Then for any $v_\theta\in\mathcal{W}$, the following estimates hold
\begin{enumerate}
\item[{\rm(i)}] $\|v_\theta\|_{L^\infty(\Omega)}\leq WB_\theta$,   $\|v_\theta-v_{\tilde{\theta}}\|_{L^\infty(\Omega)}\leq 2LW^LB_\theta^{L-1}|\theta-\tilde\theta|_{\ell^\infty}$;
\item[{\rm(ii)}] $\|\nabla v_\theta\|_{L^\infty(\Omega; \mathbb{R}^d)}\leq \sqrt{2}W^{L-1}B_\theta^L$,
$\|\nabla (v_\theta-v_{\tilde{\theta}})\|_{L^\infty(\Omega; \mathbb{R}^d)}\leq  \sqrt{2}L^2W^{2L-2}B_\theta^{2L-2}|\theta-\tilde\theta|_{\ell^\infty}$;
\item[{\rm(iii)}] $\|\Delta v_\theta\|_{L^\infty(\Omega)}\leq 2LW^{2L-2} B_\theta^{2L} $,
$\|\Delta (v_\theta-v_{\tilde{\theta}})\|_{L^\infty(\Omega)}\leq  8N_\theta L^2W^{3L-3}B_\theta^{3L-3}|\theta-\tilde\theta|_{\ell^\infty}$.
\end{enumerate}
\end{lemma}

Next we give boundedness and Lipschitz continuity of functions in $\mathcal{H}_r$, $\mathcal{H}_d$ and $\mathcal{H}_n$.

\begin{lemma}\label{lem:fcn-Lip}
There exists a constant $C=C(\|f^\delta\|_{L^\infty(\partial\Omega)},\|g^\delta\|_{L^\infty(\partial\Omega)},B_\mathbf{c},\gamma,\beta,M,d)$ such that
$ \|h( v_\theta)\|_{L^\infty(\Omega)} \le C L^2 W^{4L-4}B_\theta^{4L}$ for all $h \in \mathcal{H}_r$;  $\|h(v_\theta )\|_{L^\infty(\partial\Omega)} \le C W^2B_\theta^2$ for all $h \in \mathcal{H}_d$;
$\|h( v_\theta)\|_{L^\infty(\partial\Omega)} \le C W^{2L-2}B_\theta^{2L}$ for all $h \in \mathcal{H}_n$.
Moreover, the following Lipschitz continuity estimates hold: (i) $\|h(v_\theta) - \tilde h(v_{\tilde\theta}) \|_{L^\infty(\Omega)} \le CN_\theta L^3W^{5L-5}B_\theta^{5L-3}| \theta - \tilde\theta|_{\ell^\infty}$ for all $h,\tilde h \in \mathcal{H}_r$; (ii) $\|h(v_\theta,\mathbf{c},X )- \tilde h(v_{\tilde\theta},\tilde{\mathbf{c}},\tilde{X})\|_{L^\infty(\partial\Omega)} \le CMLW^{L+1}B_\theta^L(| \theta - \tilde \theta  |_{\ell^\infty}+| \mathbf{c} - \tilde{\mathbf{c}}  |_{\ell^\infty}+\|X - \tilde X \|_{\infty})$ for all $h,\tilde h \in \mathcal{H}_d$; (iii) $\|h(v_\theta,\mathbf{c},X)-\tilde h(v_{\tilde\theta},\tilde{\mathbf{c}},\tilde{X})\|_{L^\infty(\partial\Omega)} \le CML^2 W^{3L-3} B_\theta^{3L-2}(| \theta - \tilde \theta  |_{\ell^\infty}+| \mathbf{c} - \tilde{\mathbf{c}}  |_{\ell^\infty}+\|X - \tilde X  \|_{\infty})$ for all $h,\tilde h \in \mathcal{H}_n$.
\end{lemma}
\begin{proof}
It follows from the definition of $\Phi(\mathbf{x})$ in \eqref{fundamental-s} that
\begin{align*}
\sup_{(\mathbf{c},X)\in I_\mathbf{c}^M\times\Omega_\gamma^M}\Big\|\sum_{j=1}^Mc_j\Phi_{\mathbf{x}_j}\Big\|_{L^\infty(\partial\Omega)}&\leq c_dMB_c
\left\{
	\begin{aligned}
		\max\{|\log\gamma|,|\log\beta|\},&\quad d=2,\\
		\gamma^{2-d},&\quad d\geq3,
	\end{aligned}\right.\\
\sup_{(\mathbf{c},X)\in I_\mathbf{c}^M\times\Omega_\gamma^M}\Big\|\sum_{j=1}^Mc_j\partial_\nu\Phi_{\mathbf{x}_j}\Big\|_{L^\infty(\partial\Omega)}&\leq	 c_d'MB_c\gamma^{1-d}.
\end{align*}
Likewise, with $C=C(\beta,\gamma,B_\mathbf{c},M,d)$, we deduce
\begin{align*}
    &\Big\|\sum_{j=1}^M\left(c_j\Phi_{\mathbf{x}_j}-\tilde{c}_j\Phi_{\tilde{\mathbf{x}}_j}\right)\Big\|_{L^\infty(\partial\Omega)}\\
    \leq&\Big\|\sum_{j=1}^M\left(c_j-\tilde{c}_j\right)\Phi_{\mathbf{x}_j}\Big\|_{L^\infty(\partial\Omega)}+\Big\|\sum_{j=1}^M\tilde{c}_j(\Phi_{\mathbf{x}_j}-\Phi_{\tilde{\mathbf{x}}_j})\Big\|_{L^\infty(\partial\Omega)}\\
    \leq& C\big(| \mathbf{c} - \tilde{\mathbf{c}}  |_{\ell^\infty}+\max_{1\leq j
    \leq M}\|\mathbf{x}_j-\tilde{\mathbf{x}}_j\|\big)
    \leq C\big(| \mathbf{c} - \tilde{\mathbf{c}}  |_{\ell^\infty}+\max_{1\leq j
    \leq M}\|\mathbf{x}_j-\tilde{\mathbf{x}}_j\|_{\ell^\infty}\big).
\end{align*}
Similarly, with $C= C(\beta,\gamma,B_\mathbf{c},M,d)$,
\begin{align*}
    &\Big\|\sum_{j=1}^M\left(c_j\partial_\nu\Phi_{\mathbf{x}_j}-\tilde{c}_j\partial_\nu\Phi_{\tilde{\mathbf{x}}_j}\right)\Big\|_{L^\infty(\partial\Omega)}\\
    \leq&\Big\|\sum_{j=1}^M\left(c_j-\tilde{c}_j\right)\partial_\nu\Phi_{\mathbf{x}_j}\Big\|_{L^\infty(\partial\Omega)}+\Big\|\sum_{j=1}^M\tilde{c}_j(\partial_\nu\Phi_{\mathbf{x}_j}-\partial_\nu\Phi_{\tilde{\mathbf{x}}_j})\Big\|_{L^\infty(\partial\Omega)}\\
    \leq& C\big(| \mathbf{c} - \tilde{\mathbf{c}}  |_{\ell^\infty}+\max_{1\leq j\leq M}\|\mathbf{x}_j-\tilde{\mathbf{x}}_j\|\big)
    \leq C\big(| \mathbf{c} - \tilde{\mathbf{c}}  |_{\ell^\infty}+\max_{1\leq j\leq M}\|\mathbf{x}_j-\tilde{\mathbf{x}}_j\|_{\ell^\infty}\big).
\end{align*}
Then all the estimates are direct from Lemma \ref{lem:NN-Lip}.
\end{proof}

\bibliographystyle{abbrv}
\bibliography{refer}

\begin{thebibliography}{10}

\bibitem{anikonov2013inverse}
Y.~E. Anikonov, B.~A. Bubnov, and G.~N. Erokhin.
\newblock {\em {Inverse and Ill-posed Sources Problems}}.
\newblock Walter de Gruyter, Berlin, 2013.

\bibitem{AnthonyBartlett:1999}
M.~Anthony and P.~L. Bartlett.
\newblock {\em {Neural Network Learning: Theoretical Foundations}}.
\newblock Cambridge University Press, Cambridge, 1999.

\bibitem{BartlettMendelson:2002}
P.~L. Bartlett and S.~Mendelson.
\newblock Rademacher and {G}aussian complexities: risk bounds and structural
  results.
\newblock {\em J. Mach. Learn. Res.}, 3:463--482, 2002.

\bibitem{Berggren:2004}
M.~Berggren.
\newblock Approximations of very weak solutions to boundary-value problems.
\newblock {\em SIAM J. Numer. Anal.}, 42(2):860--877, 2004.

\bibitem{BottouCurtis:2018}
L.~Bottou, F.~E. Curtis, and J.~Nocedal.
\newblock Optimization methods for large-scale machine learning.
\newblock {\em SIAM Rev.}, 60(2):223--311, 2018.

\bibitem{ByrdLu:1995}
R.~H. Byrd, P.~Lu, J.~Nocedal, and C.~Y. Zhu.
\newblock A limited memory algorithm for bound constrained optimization.
\newblock {\em SIAM J. Sci. Comput.}, 16(5):1190--1208, 1995.

\bibitem{cai2001finite}
Z.~Cai and S.~Kim.
\newblock A finite element method using singular functions for the {P}oisson
  equation: corner singularities.
\newblock {\em SIAM J. Numer. Anal.}, 39(1):286--299, 2001.

\bibitem{CaiShin:2001}
Z.~Cai, S.~Kim, and B.-C. Shin.
\newblock Solution methods for the {P}oisson equation with corner
  singularities: numerical results.
\newblock {\em SIAM J. Sci. Comput.}, 23(2):672--682, 2001.

\bibitem{chafik2000some}
M.~Chafik, A.~El~Badia, and T.~Ha-Duong.
\newblock On some inverse {EEG} problems.
\newblock In {\em {Inverse Problems in Engineering Mechanics II}}, pages
  537--544. Elsevier, 2000.

\bibitem{CuckerSmale:2002}
F.~Cucker and S.~Smale.
\newblock On the mathematical foundations of learning.
\newblock {\em Bull. Amer. Math. Soc. (N.S.)}, 39(1):1--49, 2002.

\bibitem{Jin:DNN-Control}
Y.~Dai, B.~Jin, R.~Sau, and Z.~Zhou.
\newblock Solving elliptic optimal control problems via neural networks and
  optimality system.
\newblock {\em Preprint, arXiv:2308.11925}, 2023.

\bibitem{DuLiSun:2023}
H.~Du, Z.~Li, J.~Liu, Y.~Liu, and J.~Sun.
\newblock Divide-and-conquer {DNN} approach for the inverse point source
  problem using a few single frequency measurements.
\newblock {\em Inverse Problems}, 39(11):115006, 19, 2023.

\bibitem{ElBadiaElHajj:2012}
A.~El~Badia and A.~El~Hajj.
\newblock H\"{o}lder stability estimates for some inverse pointwise source
  problems.
\newblock {\em C. R. Math. Acad. Sci. Paris}, 350(23-24):1031--1035, 2012.

\bibitem{Badia2013}
A.~El~Badia and A.~El~Hajj.
\newblock Stability estimates for an inverse source problem of {H}elmholtz’s
  equation from single {C}auchy data at a fixed frequency.
\newblock {\em Inverse Problems}, 29(12):125008, 20, 2013.

\bibitem{el1998some}
A.~El~Badia and T.~Ha-Duong.
\newblock Some remarks on the problem of source identification from boundary
  measurements.
\newblock {\em Inverse Problems}, 14(4):883--891, 1998.

\bibitem{el2000inverse}
A.~El~Badia and T.~Ha-Duong.
\newblock An inverse source problem in potential analysis.
\newblock {\em Inverse Problems}, 16(3):651--663, 2000.

\bibitem{evans10}
L.~C. Evans.
\newblock {\em {Partial Differential Equations}}.
\newblock AMS, Providence, RI, second edition, 2010.

\bibitem{Fix:1973}
G.~J. Fix, S.~Gulati, and G.~I. Wakoff.
\newblock On the use of singular functions with finite element approximations.
\newblock {\em J. Comput. Phys.}, 13:209--228, 1973.

\bibitem{Fries:2010}
T.-P. Fries and T.~Belytschko.
\newblock The extended/generalized finite element method: an overview of the
  method and its applications.
\newblock {\em Internat. J. Numer. Methods Engrg.}, 84(3):253--304, 2010.

\bibitem{Gautschi1962OnIO}
W.~Gautschi.
\newblock On inverses of {V}andermonde and confluent {V}andermonde matrices.
\newblock {\em Numer. Math.}, 4:117--123, 1962.

\bibitem{gjerde2019splitting}
I.~G. Gjerde, K.~Kumar, J.~M. Nordbotten, and B.~Wohlmuth.
\newblock Splitting method for elliptic equations with line sources.
\newblock {\em ESAIM: Math. Model. Numer. Anal.}, 53(5):1715--1739, 2019.

\bibitem{GruterWidman:1982}
M.~Gr\"{u}ter and K.-O. Widman.
\newblock The {G}reen function for uniformly elliptic equations.
\newblock {\em Manuscripta Math.}, 37(3):303--342, 1982.

\bibitem{GuhringRaslan:2021}
I.~Guhring and M.~Raslan.
\newblock Approximation rates for neural networks with encodable weights in
  smoothness spaces.
\newblock {\em Neural Networks}, 134:107--130, 2021.

\bibitem{HuJinZhou:2022}
T.~Hu, B.~Jin, and Z.~Zhou.
\newblock Solving elliptic problems with singular sources using singularity
  splitting deep {R}itz method.
\newblock {\em SIAM J. Sci. Comput.}, 45(4):A2043--A2074, 2023.

\bibitem{hu2023solving}
T.~Hu, B.~Jin, and Z.~Zhou.
\newblock Solving {P}oisson problems in polygonal domains with singularity
  enriched physics informed neural networks.
\newblock {\em SIAM J. Sci. Comput.}, 46(4):C369--C398, 2024.

\bibitem{ikehata1999reconstruction}
M.~Ikehata.
\newblock Reconstruction of obstacles from boundary measurements.
\newblock {\em Wave Motion}, 30(3):205--223, 1999.

\bibitem{isakov1990inverse}
V.~Isakov.
\newblock {\em {Inverse Source Problems}}.
\newblock AMS, Providence, RI, 1990.

\bibitem{JiaoLai:2022cicp}
Y.~Jiao, Y.~Lai, D.~Li, X.~Lu, F.~Wang, Y.~Wang, and J.~Z. Yang.
\newblock A rate of convergence of physics informed neural networks for the
  linear second order elliptic {PDE}s.
\newblock {\em Commun. Comput. Phys.}, 31(4):1272--1295, 2022.

\bibitem{JinLiLu:2022}
B.~Jin, X.~Li, and X.~Lu.
\newblock Imaging conductivity from current density magnitude using neural
  networks.
\newblock {\em Inverse Problems}, 38(7):075003, 36, 2022.

\bibitem{Karniadakis:2021nature}
G.~E. Karniadakis, I.~G. Kevrekidis, L.~Lu, P.~Perdikaris, S.~Wang, and
  L.~Yang.
\newblock Physics-informed machine learning.
\newblock {\em Nature Rev. Phys.}, 3:422--440, 2021.

\bibitem{KingmaBa:2015}
D.~P. Kingma and J.~Ba.
\newblock Adam: A method for stochastic optimization.
\newblock In {\em 3rd International Conference for Learning Representations},
  San Diego, 2015.

\bibitem{kriegsmann1988source}
G.~Kriegsmann and W.~Olmstead.
\newblock Source identification for the heat equation.
\newblock {\em Appl. Math. Lett.}, 1(3):241--245, 1988.

\bibitem{LiLiangWang:2023}
P.~Li, Y.~Liang, and Y.~Wang.
\newblock A data-assisted two-stage method for the inverse random source
  problem.
\newblock {\em SIAM J. Imaging Sci.}, 16(4):1929--1952, 2023.

\bibitem{LuChenLu:2021}
Y.~Lu, H.~Chen, J.~Lu, L.~Ying, and J.~Blanchet.
\newblock Machine learning for elliptic {PDE}s: Fast rate generalization bound,
  neural scaling law and minimax optimality.
\newblock In {\em International Conference on Learning Representations}, 2022.

\bibitem{lu2021priori}
Y.~Lu, J.~Lu, and M.~Wang.
\newblock A priori generalization analysis of the deep {R}itz method for
  solving high dimensional elliptic partial differential equations.
\newblock In {\em Conference on Learning Theory}, pages 3196--3241. PMLR, 2021.

\bibitem{michel2004128}
C.~M. Michel, G.~Lantz, L.~Spinelli, R.~G. De~Peralta, T.~Landis, and M.~Seeck.
\newblock 128-channel {EEG} source imaging in epilepsy: clinical yield and
  localization precision.
\newblock {\em J. Clinical Neurophys.}, 21(2):71--83, 2004.

\bibitem{MICHEL20042195}
C.~M. Michel, M.~M. Murray, G.~Lantz, S.~Gonzalez, L.~Spinelli, and R.~{Grave
  de Peralta}.
\newblock {EEG} source imaging.
\newblock {\em Clinical Neurophys.}, 115(10):2195--2222, 2004.

\bibitem{Mohri:2018}
M.~Mohri, A.~Rostamizadeh, and A.~Talwalkar.
\newblock {\em {Foundations of Machine Learning}}.
\newblock MIT Press, Cambridge, MA, 2018.

\bibitem{Rado:1949}
R.~Rado.
\newblock Factorization of even graphs.
\newblock {\em Quart. J. Math.}, os-20(1):95--104, 1949.

\bibitem{sirignano2018dgm}
J.~Sirignano and K.~Spiliopoulos.
\newblock Dgm: A deep learning algorithm for solving partial differential
  equations.
\newblock {\em J. Comput. Phys.}, 375:1339--1364, 2018.

\bibitem{WangPerdikaris:2022jcp}
S.~Wang, X.~Yu, and P.~Perdikaris.
\newblock When and why {PINN}s fail to train: a neural tangent kernel
  perspective.
\newblock {\em J. Comput. Phys.}, 449:110768, 28, 2022.

\bibitem{wolf2018mathematical}
M.~M. Wolf.
\newblock {Mathematical Foundations of Supervised Learning}.
\newblock Lecture notes, available at
  \url{https://www-m5.ma.tum.de/foswiki/pub/M5/Allgemeines/MA4801_2021S/ML.pdf}
  (accessed on Feb. 1, 2023), 2021.

\bibitem{Zhang2023}
H.~Zhang and J.~Liu.
\newblock Solving an inverse source problem by deep neural network method with
  convergence and error analysis.
\newblock {\em Inverse Problems}, 39(7):075013, 2023.

\bibitem{ZhangLiLiu:2023}
M.~Zhang, Q.~Li, and J.~Liu.
\newblock On stability and regularization for data-driven solution of parabolic
  inverse source problems.
\newblock {\em J. Comput. Phys.}, 474:111769, 20, 2023.

\end{thebibliography}

\end{document}